\documentclass[12pt]{amsart}
\usepackage{amsmath,amsthm,amsfonts,amssymb,mathrsfs}
\usepackage[all]{xy}
\usepackage{eufrak}
\usepackage{amssymb}
\usepackage{latexsym}
\usepackage{amsmath,amsthm}
\usepackage{aeguill}
\usepackage{amssymb}
\usepackage{mathrsfs}
\usepackage{hyperref}
\usepackage{float}
\usepackage{graphicx}
\usepackage{etoolbox}
\usepackage{enumerate}
\usepackage{ textcomp }
\usepackage[total={6.5in,8.75in}, top=1.2in, left=1.0in, includefoot]{geometry}
\usepackage[utf8]{inputenc}
\usepackage{tikz}[2013/12/13]

\makeatletter
\pretocmd{\chapter}{\addtocontents{toc}{\protect\addvspace{15\p@}}}{}{}
\pretocmd{\section}{\addtocontents{toc}{\protect\addvspace{3\p@}}}{}{}
\makeatother

\makeatletter
\def\@tocline#1#2#3#4#5#6#7{\relax
  \ifnum #1>\c@tocdepth 
  \else
    \par \addpenalty\@secpenalty\addvspace{#2}%
    \begingroup \hyphenpenalty\@M
    \@ifempty{#4}{%
      \@tempdima\csname r@tocindent\number#1\endcsname\relax
    }{%
      \@tempdima#4\relax
    }%
    \parindent\z@ \leftskip#3\relax \advance\leftskip\@tempdima\relax
    \rightskip\@pnumwidth plus4em \parfillskip-\@pnumwidth
    #5\leavevmode\hskip-\@tempdima
      \ifcase #1
       \or\or \hskip .5em \or \hskip 1em \else \hskip 1.5em \fi%
      #6\nobreak\relax
    \dotfill\hbox to\@pnumwidth{\@tocpagenum{#7}}\par
    \nobreak
    \endgroup
  \fi}
\makeatother
\newcommand{\C}{\mathbb{C}}
\newcommand{\N}{\mathbb{N}}
\newcommand{\Z}{\mathbb{Z}}
\newcommand{\R}{\mathbb{R}}
\newcommand{\Q}{\mathbb{Q}}

\newcommand{\A}{\mathbb{A}}
\newcommand{\X}{\mathbb{X}}

\newcommand{\Gal}{\operatorname{Gal}}

\newcommand{\Frob}{\mathrm{Frob}}

\newcommand{\ab}{\operatorname{ab}}

\renewcommand{\ss}{\operatorname{ss}}
\newcommand{\der}{\operatorname{der}}

\newcommand{\End}{\operatorname{End}}

\newcommand{\hto}{\hookrightarrow}

\newcommand{\GL}{\mathrm{GL}}

\newcommand{\SL}{\mathrm{SL}}

\newcommand{\Char}{\mathrm{Char}}

\newcommand{\bG}{\mathbf{G}}
\newcommand{\bH}{\mathbf{H}}

\newcommand{\bN}{\mathbf{N}}
\newcommand{\bS}{\mathbf{S}}
\newcommand{\bT}{\mathbf{T}}

\newcommand{\bY}{\mathbf{Y}}

\newtheorem{thm}{Theorem}[section]

\newtheorem{cor}[thm]{Corollary}

\newtheorem{prop}[thm]{Proposition}

\newtheorem{remark}[thm]{Remark}

\theoremstyle{definition}
\newtheorem{definition}[thm]{Definition}
  \newtheorem{example}[thm]{Example}

\newtheorem{innercustomconj}{Conjecture}
\newenvironment{customconj}[1]
  {\renewcommand\theinnercustomconj{#1}\innercustomconj}
  {\endinnercustomconj}
\begin{document}

\title[The abelian part of a compatible system]{The abelian part of a compatible system and $\ell$-independence of the Tate conjecture}

\author{Chun Yin Hui}

\date{August, 2018}
\thanks{The present project was supported partly by the National Research Fund, Luxembourg, cofunded under the Marie Curie Actions of the European Commission (FP7-COFUND) and partly by China's Thousand Talents Plan: The Recruitment Program for Young Professionals.}

\maketitle

\begin{abstract} 
Let $K$ be a number field and $\{V_\ell\}_\ell$ be a rational strictly compatible system 
of semisimple Galois representations of $K$ arising from geometry.
Let $\bG_\ell$ and $V_\ell^{\ab}$ be respectively the algebraic monodromy group
and the maximal abelian subrepresentation of $V_\ell$ for all $\ell$. 
We prove that the system $\{V_\ell^{\ab}\}_\ell$ is also a rational strictly compatible system
under some group theoretic conditions, 
e.g., when $\bG_{\ell'}$ is connected and satisfies \emph{Hypothesis A} for some prime $\ell'$.
As an application, we prove that the Tate conjecture for abelian variety $X/K$ is independent of $\ell$ 
if the algebraic monodromy groups of the Galois representations of $X$
satisfy the required conditions.
\end{abstract}

\newpage
\section{Introduction}
\subsection{Conjectures}
Fix a number field $K$ and an embedding $\overline{\Q}_{\ell}\hookrightarrow\C$ for every prime $\ell$.
An $\overline\Q_\ell$-representation $\overline W_\ell$ of $\Gal(\overline K/K)$ is said to be \emph{arising from geometry}
if it is isomorphic to a subquotient of $\bigoplus_{j=1}^m H^{i_j}(\overline{X}_{j},\overline{\Q}_{\ell}(r_j))$, 
where $X_j$ is smooth projective over $K$, $\overline{X}_{j}:=X_j\times\overline K$, $m\in\N$, $1\leq j\leq m$, and $i_j,r_j\in \Z$.
The Langlands reciprocity conjecture \cite[$\mathsection2$]{Lan79} 
(see also \cite[$\mathsection4$]{Cl90},\cite[Conj. 3.5, 3.4]{Tay04})
predicts that if $\overline{W}_{\ell'}$ is an  
$n$-dimensional irreducible $\overline{\Q}_{\ell'}$-representation of $\Gal(\overline K/K)$ arising from geometry,
then there is a cuspidal automorphic representations $\pi$ of $\GL_n(\A_K)$
preserving the local L-factors at almost all primes $p$, i.e., 
$L_p(\pi,p^{-s})=L_p(\overline{W}_{\ell'}\otimes\C, p^{-s})$.
Moreover, for each prime $\ell$, one associates to $\pi$ a $\overline\Q_\ell$-representation $\overline W_\ell$ of $\Gal(\overline K/K)$
such that the system $\{\overline W_\ell\}_\ell$ is weakly compatible\footnote{The weakly compatibility here
is much weaker than the weakly compatibility in \cite[$\mathsection1$]{Tay04}.} 
in the following sense.


\begin{definition}\label{wcompat}
Let $\{\rho_\iota:\Gal(\overline K/K)\to \GL_n(F_\iota)\}_\iota$ be a family of (continuous) Galois representations,
where  $F_\iota$ is an algebraic extension of some
characteristic zero non-Archimedean local field and
$\iota$ belongs to an index set $I$. The system $\{\rho_\iota\}_\iota$ is said to be \emph{weakly compatible} if there exists an embedding
$F_\iota\hto\C$ for every $\iota$ such that:
\begin{enumerate}[(i)]
\item for each $\iota$, $\rho_\iota$ is unramified outside a finite set $\Sigma_\iota$ of finite places $v$ of $K$; 
\item for every pair $\iota,\iota'$ of indices, 
the characteristic polynomials $P_{v,\rho_\iota}(t),P_{v,\rho_{\iota'}}(t)$ of the Frobenii above $v$ agree in $\C[t]$ 
 via $F_\iota\hto\C \hookleftarrow F_{\iota'}$ if $v$ a finite place of $K$ outside $\Sigma_\iota\cup \Sigma_{\iota'}$.
\end{enumerate} 
\end{definition}

The following more refined notion of compatibility of Galois representations (Definition \ref{compat}) was 
introduced by Serre in \cite[Ch. I $\mathsection2$]{Se98}.

\begin{definition}\label{rat}
Let $E_\ell$ be a finite extension of $\Q_\ell$, $E$ be a number field contained in $E_\ell$, and
$\rho_\ell:\Gal(\overline K/K)\to \GL_n(E_\ell)$ be a continuous $E_\ell$-representation.
Then $\rho_\ell$ is said to be \emph{$E$-rational} if $\rho_\ell$ is unramified
outside a finite subset $S$ of finite places $v$ of $K$ and the characteristic polynomial $P_{v,\rho_\ell}(t)$ of the Frobenii above $v$
belongs to $E[t]$ for all $v$ outside $S$. The term $\Q$-rational is abbreviated as \emph{rational}.
\end{definition}

For any finite place $v$ of a number field, let $p_v$ be the characteristic of $v$.
For any number field $E$ with a finite place $\lambda$, let $E_\lambda$ be the completion of $E$ with respect to $\lambda$.
An $E_\lambda$-representation is also called a $\lambda$-adic representation.

\begin{definition}\label{compat}
For each finite place $\lambda$ of $E$, let $\rho_\lambda$ be an $E$-rational $\lambda$-adic
representation of $\Gal(\overline K/K)$. The system $\{\rho_\lambda\}_\lambda$ is said to be \emph{strictly compatible} if there exists 
a finite subset $S$ of finite places $v$ of $K$ such that: 
\begin{enumerate}[(i)]
\item let $S_\lambda:=\{v:~p_v=p_\lambda\}$, then $\rho_\lambda$ is unramified if $v\notin S\cup S_\lambda$.
\item $P_{v,\rho_\lambda}(t)=P_{v,\rho_{\lambda'}}(t)$ if $v\notin S\cup S_\lambda\cup S_{\lambda'}$.
\end{enumerate}
The smallest possible $S$ is called the \emph{exceptional set} of the system.
\end{definition}

A system of Galois representations is said to be \emph{semisimple} if every member is a semisimple representation. 
 Let $X$ be smooth projective over $K$. 
The Galois representation $H^i(\overline{X},\Q_\ell)$ 
is conjectured by Grothendieck-Serre to be semisimple for all $\ell$ \cite[p. 109]{Ta65}.
Deligne has proven that the $\ell$-adic representations $\{V_\ell=H^i(\overline{X},\Q_\ell)\}_\ell$ form a rational 
strictly compatible system with exceptional set equal to 
the finite places $v$ of $K$ for which $X$ does not have good reduction \cite{De74}.
The semisimplicity conjecture, the conjectural correspondence of Langlands, and the 
compatibility condition lead us to the following conjecture on the decomposition of the
rational strictly compatible system $\{V_\ell\}_\ell$ arising from geometry, i.e., 
 $\{V_\ell\}_\ell$ is a rational strictly compatible system 
and $V_\ell\otimes\overline\Q_\ell$ arises from geometry for all primes $\ell$.

\begin{customconj}{De}\label{D}
\textit{Let $\{V_\ell\}_\ell$ be a rational strictly compatible system arising from geometry.
Then there exist weakly compatible systems of irreducible $\overline{\Q}_\ell$-representations $\{\overline{W}^{(k)}_\ell\}_\ell$ for $1\leq k\leq m$ 
such that for all $\ell$, $V_\ell$ admits the following decomposition:
$$V_\ell\otimes\overline{\Q}_\ell\cong \overline{W}^{(1)}_\ell\oplus \overline{W}^{(2)}_\ell\oplus \cdots \oplus \overline{W}^{(m)}_\ell.$$}
\end{customconj}

Since both the semisimplicity conjecture and the Langlands reciprocity conjecture are wide open\footnote{ 
The function field analogue of the Langlands conjecture (for $\GL_n$) is proved by Lafforgue \cite{La02}.
For recent advances toward the Langlands conjecture, see \cite{BLGGT14}.}, 
Conjecture \ref{D} in general is out of reach, 
especially when the irreducibles of $V_\ell\otimes\overline{\Q}_\ell$ 
are not all of dimension one or two. 
However, the next conjecture, which is a special case of the above one, 
can be established via $\ell$-independence of $\{V_\ell\}_\ell$ in some cases.

\begin{definition}\label{abpart}
Let $V$ be a finite dimensional semisimple representation over a field of a group $G$. 
The \emph{abelian part} $V^{\mathrm{ab}}$ of $V$ is defined 
to be the maximal abelian subrepresentation of  $V$.
\end{definition}

\begin{customconj}{Ab}\label{A}
\textit{Let $\{V_\ell\}_\ell$ be a semisimple 
rational strictly compatible system arising from geometry, 
with exceptional set $S$.
Then the abelian part $\{V_\ell^{\mathrm{ab}}\}_\ell$ of $\{V_\ell\}_\ell$ is also 
rational strictly compatible with exceptional set contained in $S$.}
\end{customconj}

\subsection{Theorems} 
We describe the main results of this paper.
Let $\{V_\ell\}_\ell$ be a semisimple rational strictly compatible system 
of $\ell$-adic representations of $\Gal(\overline{K}/K)$ arising from geometry ($\mathsection1.1$). 
Suppose the representation $V_\ell$ is represented by 
$$\varphi_\ell:\Gal(\overline{K}/K)\to \GL(V_\ell)\cong\GL_N(\Q_\ell)$$
for all $\ell$. The \emph{algebraic monodromy group} of $V_\ell$, denoted by $\bG_\ell$, is defined as the Zariski closure of the Galois image $\varphi_\ell(\Gal(\overline{K}/K))$ in the $\Q_\ell$-algebraic group $\GL_{V_\ell}$. Since $\varphi_\ell$ is semisimple, 
we identify $\bG_\ell$ as a reductive subgroup of $\GL_{N,\Q_\ell}$ for all $\ell$.
The tautological representation $\bG_\ell^\circ\hto \GL_{N,\Q_\ell}$ is conjectured to be independent of $\ell$ \cite{Se94}.
For results in this direction, see 
\cite{Se72,Se85},\cite{Ri76},\cite{Ch92},\cite{Pi98},\cite{LP95} for abelian varieties and 
\cite{LP92},\cite{Chi04},\cite{Hu13,Hu15,Hu18} for more general cases.

Denote by $\bG_\ell^{\der}$
the derived group  of the identity component of $\bG_\ell$. 
Then $\bG_\ell^{\der}$ is connected and semisimple.
Embed $\Q_\ell$ into $\C$ for all $\ell$ so that the complex semisimple Lie groups $\bG_\ell^{\der}(\C)$ are all identified
with subgroups of $\GL_N(\C)$.
Let $\bT_{\ell,\C}^{\ss}$ be a maximal torus of $\bG_\ell^{\der}(\C)$ for some $\ell$ 
and $\X_\ell^{\ss}$ be the character group
of $\bT_{\ell,\C}^{\ss}$. Then the automorphism group of $\bT_{\ell,\C}^{\ss}$ is canonically isomorphic 
to the automorphism group of $\X^{\ss}_\ell$.
Let $R_\ell\subset\X_\ell^{\ss}$ be the set of \emph{roots} of $\bG_\ell^{\der}(\C)$ with respect to $\bT_{\ell,\C}^{\ss}$.
The motivation of the paper is to investigate Conjecture \ref{A} by properties of $\bG_\ell$:
\vspace{.05in}

\begin{center}
\begin{tabular}{ccc} 
Connectedness of $\bG_\ell$ &&\\
+ & & \\
$\ell$-independence of $\bG_\ell^{\der}(\C)\hto\GL_N(\C)$ & $\Rightarrow$ & Conjecture \ref{A} for $\{V_\ell\}_\ell$. \\
+ &  & \\
an invariance of roots $R_\ell$ criterion &&
\end{tabular}
\end{center}

\begin{thm}\label{thm1}
Let $\{V_\ell\}_\ell$ be a semisimple rational strictly compatible system arising from geometry, with exceptional set $S$. 
Let $\bG_\ell$ be the algebraic monodromy group of $V_\ell$ for all $\ell$.
Suppose the following conditions hold:
\begin{enumerate}[(i)]
\item $\bG_\ell$ is connected for all $\ell$;
\item the conjugacy class of $\bG_\ell^{\der}(\C)$ in $\GL_N(\C)$ is independent of $\ell$;
\item the set of roots $R_\ell$ of $\bG_\ell^{\der}(\C)$ is stable under the normalizer $N_{\GL_N(\C)}(\bT_{\ell,\C}^{\ss})$.
\end{enumerate}
Then the abelian part $\{V_\ell^{\ab}\}_\ell$ is rational strictly compatible with exceptional set contained in $S$.
\end{thm}

The connectedness condition \ref{thm1}(i) is independent of $\ell$ by Serre (see $\mathsection2.1$). 
Thus, one can take a finite extension of $K$ for \ref{thm1}(i) to hold. 
The invariance of roots criterion \ref{thm1}(iii) is independent of the choice of the maximal torus $\bT_{\ell,\C}^{\ss}$
of $\bG_\ell^{\der}(\C)$.
The $\ell$-independence condition \ref{thm1}(ii) implies that the criterion \ref{thm1}(iii) is independent of $\ell$.
The condition \ref{thm1}(ii) is a consequence of a more general (Mumford-Tate) conjecture asserting that the $\Q_\ell$-representations
$\bG_\ell^\circ\hookrightarrow\GL_{V_\ell}$ admit a common $\Q$-model for all $\ell$.
We provide explicit group theoretic conditions for Theorem \ref{thm1}.

\begin{thm}\label{thm2}
Let $\{V_\ell\}_\ell$ be a semisimple rational strictly compatible system arising from geometry, with exceptional set  $S$. 
Let $\bG_\ell$ be the algebraic monodromy group of $V_\ell$ for all $\ell$.
Suppose (i) $\bG_\ell$ is connected for all $\ell$ and one of the following conditions hold:
\begin{enumerate}
\item[(ii)] Hypothesis A: there exists some $\ell'$ such that 
$\bG_{\ell'}^{\der}(\C)$ has at most one factor of type $A_4$ and 
if $\bH_\C$ is an almost simple factor of $\bG_{\ell'}^{\der}(\C)$, then $\bH_\C$ is of type $A_n$ for some
$n\in\N\backslash\{1,2,3,5,7,8\}$;
\item[(ii')] the conjugacy class of $\bG_\ell^{\der}(\C)$ in $\GL_N(\C)$ is independent of $\ell$ 
and $\bG_\ell^{\der}(\C)$ is almost simple of type different from $A_7,A_8,B_4,D_8$.
\end{enumerate}
Then the conditions of Theorem \ref{thm1} hold.
\end{thm}

Under Hypothesis A, one can prove the following for Conjecture \ref{D}.

\begin{thm}\label{thm3}
Let $\{V_\ell\}_\ell$ be a semisimple rational strictly compatible system arising from geometry. 
Let $\bG_\ell$ be the algebraic monodromy group of $V_\ell$ for all $\ell$.
Suppose (i) $\bG_\ell$ is connected for all $\ell$ and (ii) Hypothesis A is satisfied.
Let $\overline{W}_{\ell'}^{\oplus m}$ be a $\overline{\Q}_{\ell'}$-subrepresentation of $V_{\ell'}\otimes\overline{\Q}_{\ell'}$ 
for some prime $\ell'$, some $m\in\N$, and
$\{\overline{W}_\ell\}_\ell$ be a 
weakly compatible system of irreducible $\overline{\Q}_\ell$-representations satisfying the following conditions:
\begin{enumerate}
\item[(iii)] $\overline{W}_{\ell'}$ is a member of the system $\{\overline{W}_\ell\}_\ell$;
\item[(iv)] for all $\ell$, descend $\overline{W}_\ell$ to an $E_\ell$-representation $W_\ell$ in which $E_\ell$ is 
a finite extension of $\Q_\ell$. Then the local representations of $W_\ell$ above $\ell$ are of Hodge-Tate type (Definition \ref{HT}).
\end{enumerate}
Then $\{\overline{W}_\ell^{\oplus m}\}_\ell$ is a subsystem of $\{V_\ell\otimes\overline{\Q}_\ell\}_\ell$.
\end{thm}

\subsection{Applications}
Let $X$ be a smooth projective variety defined over $K$ and 
$V_\ell:=H^{2r}(\overline{X},\Q_\ell(r))$ be 
the $\ell$-adic representation of $\Gal(\overline K/K)$.
Denote by $\mathfrak{g}_\ell$ the Lie algebra of the image of $\Gal(\overline K/K)$ in $\GL(V_\ell)$ and by $CH^r(\overline{X})$
the Chow group of codimension $r$ algebraic cycles of $\overline{X}$. 
The Tate conjecture \cite{Ta65} asserts that 
the image of the cycle class map
$$c_\ell\otimes1: CH^r(\overline{X})\otimes_\Z\Q_\ell\to H^{2r}(\overline{X},\Q_\ell(r))$$
is equal to the space of invariants $H^{2r}(\overline{X},\Q_\ell(r))^{\mathfrak{g}_\ell}$. 
The Tate conjecture is also wide open and not known to be independent of $\ell$.

\begin{cor}\label{indepTate} (to Conjecture \ref{A})
Assume the Galois representation $V_\ell:=H^{2r}(\overline{X},\Q_\ell(r))$ 
is semisimple for all $\ell$. Then Conjecture \ref{A} implies that 
the Tate conjecture for codimension $r$ algebraic cycles of $\overline X$ is independent of $\ell$.
\end{cor}

\begin{cor}\label{TateAV} (to Theorems \ref{thm1} and \ref{thm2})
Let $X/K$ be an abelian variety and $\bG_\ell$ be the algebraic monodromy group of $H^{1}(\overline{X},\Q_\ell)$.
If \ref{thm2}(ii) or \ref{thm2}(ii') holds for the groups $\bG_\ell$, 
then the Tate conjecture for codimension $r$ algebraic cycles of $\overline X$ is independent of $\ell$ for any $r\geq0$.
\end{cor}

\begin{example} Let $X_i/K$ be an absolutely simple abelian variety of prime dimension $p_i\geq 7$ such that the
endomorphism algebra $\End(\overline X_i)\otimes\Q$ is an imaginary quadratic field, $1\leq i\leq m$.
If $\bG_{i,\ell}$ denotes 
the algebraic monodromy group of the representation $H^{1}(\overline{X}_i,\Q_\ell)$, 
then the Lie type of $\bG_{i,\ell}^{\der}(\C)$ is $A_{p_i-1}$ by \cite{Ch91}, $1\leq i\leq m$.
Let $X/K$ be an abelian variety isogenous to $X_1\times X_2\times\cdots\times X_m$ and $\bG_\ell$ be 
the algebraic monodromy group of $H^{1}(\overline{X},\Q_\ell)$.
Then the Lie type of the semisimple $\bG_\ell^{\der}(\C)$ is a direct factor of 
$A_{p_1-1}\times A_{p_2-1}\times\cdots\times A_{p_m-1}$, which satisfies \ref{thm2}(ii).
\end{example}

\subsection{Idea and structure of the paper} 
The idea behind Theorems \ref{thm1} and \ref{thm3} 
is that an irreducible representation $\overline{W}_{\ell'}$ 
is a subrepresentation of a semisimple $\overline{V}_{\ell'}$ (denoted by $\overline{W}_{\ell'}\leq \overline{V}_{\ell'}$)
if and only if $\overline{V}_{\ell'}\otimes \overline{W}_{\ell'}^*=\mathrm{Hom}(\overline{W}_{\ell'},\overline{V}_{\ell'})$ has non-trivial invariants; to 
do the same thing for a family $\{\overline{W}_\ell,\overline{V}_\ell\}_\ell$ assuming $\overline{W}_{\ell'}\leq \overline{V}_{\ell'}$ for some $\ell'$, 
it suffices to show that
the tautological representation  of the algebraic monodromy group $\bH_\ell$ of $\overline{V}_\ell\otimes \overline{W}_\ell^*$ is independent of $\ell$.
Section $2$ is devoted to the proofs of the results in this section, which rely 
on some results of Serre on Galois representations and the techniques we have developed in \cite{Hu13,Hu18}. 
The appendix studies the criterion on the roots in \ref{thm1}(iii)
by a geometric configuration of the root system.
The results will complete the proof of Theorem \ref{thm2}.

\section{The abelian part of a compatible system}

\subsection{Frobenius torus} Let $E_\ell$ be a finite extension of $\Q_\ell$ and fix an
embedding $E_\ell\hto\C$ for all $\ell$. 
Let $\{V_\ell\}_\ell$ be a weakly compatible (Definition \ref{wcompat}) 
system of semisimple $E_\ell$-representations:
\begin{equation}
\psi_\ell:\Gal(\overline K/K)\to \GL( V_\ell),\hspace{.1in}\forall \ell.
\end{equation}
The algebraic monodromy group $\bG_\ell$ of $V_\ell$ is reductive. 
The method of \emph{Frobenius torus}, pioneered by Serre, is foundational 
to the study of $\ell$-independence of $\bG_\ell$. Denote by $v$ a finite place
of $K$ and $\overline v$ a place of $\overline K$ dividing $v$. Let $\psi_\ell(\Frob_{\overline v})$ be
the image of the Frobenius at $\overline v$ whenever $\psi_\ell$ is unramified at $v$.

\begin{definition}
The \emph{Frobenius torus} $\bT_{\overline v,\ell}$ is defined as the identity component of the smallest
algebraic subgroup $\bS_{\overline v,\ell}$ of $\bG_\ell$ containing the semisimple part $\psi_\ell(\Frob_{\overline v})_{\ss}$ 
of $\psi_\ell(\Frob_{\overline v})$.
\end{definition}

\begin{thm}\label{Frob} (Serre \cite{Se81}) (see also \cite[Theorem 3.7]{Ch92},\cite[Theorem 1.2]{LP97}) 
Fix a prime $\ell'$ and suppose $\bG_{\ell'}$ is connected. Denote by $p_v$ the characteristic of $v$ and by 
$q_v$ the cardinality of the residue field of $v$. 
Suppose the eigenvalues $\alpha$ of $\psi_{\ell'}(\mathrm{Frob}_{\overline v})$ are algebraic for almost all $v$ 
and for each eigenvalue $\alpha$ the following conditions hold:
\begin{enumerate}[(i)]
\item the absolute values of $\alpha$ in all complex embeddings are equal;
\item $\alpha$ is a unit at any finite place not above $p_v$;
\item for any non-Archimedean valuation $w$ of $\overline \Q$ such that $w(p_v)>0$, the ratio $w(\alpha)/w(q_v)$
belongs to a finite subset of $\Q$ that is independent of $\overline v$.
\end{enumerate}
Then there exists a proper closed subvariety $\bY$ of $\bG_{\ell'}$ such that $\bT_{\overline v,\ell'}$
is a maximal torus of $\bG_{\ell'}$ whenever $\psi_{\ell'}(\mathrm{Frob}_{\overline v})\in\bG_{\ell'}\backslash \bY$.
\end{thm}

\begin{remark}\label{FCremark}
If $\bT_{\overline v,\ell'}$ is a maximal torus of the connected reductive $\bG_{\ell'}$, 
then the centralizer of $\bT_{\overline v,\ell'}$ in $\bG_{\ell'}$ is itself.
Thus, $\bS_{\overline v,\ell'}=\bT_{\overline v,\ell'}$. Moreover,  
as the unipotent part of $\psi_{\ell'}(\mathrm{Frob}_{\overline v})$ commutes with $\bT_{\overline v,\ell'}$,
the element $\psi_{\ell'}(\mathrm{Frob}_{\overline v})$ is semisimple.
\end{remark}

The following consequences are well-known.

\begin{cor}\label{connect}
Suppose the hypotheses of Theorem \ref{Frob} hold. Then $\bG_\ell$ is connected for all $\ell$.
\end{cor}

\begin{proof}
By the Chebotarev density theorem, $\bG_\ell$ is the Zariski closure of the image of the Frobenii $\psi_{\ell}(\mathrm{Frob}_{\overline v})$
for a Dirichlet density $1$ subset $S_1$ of places $v$. Moreover, we can shrink $S_1$ such that $\bT_{\overline v,\ell'}$
is a maximal torus of $\bG_{\ell'}$ for all $v\in S_1$ by Theorem \ref{Frob}. Since $\bG_{\ell'}$ is connected reductive, 
Remark \ref{FCremark} implies that $\bS_{\overline v,\ell'}=\bT_{\overline v,\ell'}$ is connected if $v\in S_1$.
Hence, $\bG_\ell$ is generated by the connected subgroups $\bH_{\overline v,\ell}$ (for all $v\in S_1$) by the weak compatibility, where $\bH_{\overline v,\ell}$ is the Zariski closure of $\psi_{\ell}(\mathrm{Frob}_{\overline v})$ in $\bG_\ell$. We conclude that $\bG_\ell$ is connected for all $\ell$.
\end{proof}

\begin{cor}\label{FC}
Suppose the hypotheses of Theorem \ref{Frob} hold. Let $\bT_{\ell,\C}$ be a maximal torus of $\bG_\ell(\C)$. Then the tautological representation $\bT_{\ell,\C}\hookrightarrow\GL(V_\ell\otimes\C)$ is independent of $\ell$. In particular, the rank of $\bG_\ell$ is independent of $\ell$.
\end{cor}

\begin{proof}
For any $\ell,\ell'$, we can find some $\overline v$ such that both $\bT_{\overline v,\ell}$ and $\bT_{\overline v,\ell'}$ 
are maximal Frobenius tori of $\bG_\ell$ and $\bG_{\ell'}$ respectively by Corollary \ref{connect} and Theorem \ref{Frob}.
Since $\psi_{\ell}(\mathrm{Frob}_{\overline v})$ and $\psi_{\ell'}(\mathrm{Frob}_{\overline v})$ 
have the same characteristic polynomial by compatibility and are both semisimple by Remark \ref{FCremark}, 
the tautological representations of $\bT_{\overline v,\ell}(\C)$ and $\bT_{\overline v,\ell'}(\C)$ are isomorphic,
i.e., there exists an $\C$-vector space isomorphism $V_\ell\otimes\C\cong V_{\ell'}\otimes\C$ inducing the commutative diagram:
\begin{equation*}
\xymatrix{ \bT_{\overline v,\ell}(\C)   \ar@{^{(}->}[r] \ar[d]^\cong & \GL(V_\ell\otimes\C)  \ar[d]^\cong \\
            \bT_{\overline v,\ell'}(\C) \ar@{^{(}->}[r] & \GL(V_{\ell'}\otimes\C).}
						\end{equation*}
Since all maximal tori of a connected reductive group (over an algebraically closed field) are conjugate, we are done. 
\end{proof}

\begin{prop}\label{eigenvalue} \cite{De80},\cite{Se81} (see also \cite{LP97})
The conditions on Frobenius eigenvalues of Theorem \ref{Frob} 
always hold for $E_\ell$-representations arising from geometry ($\mathsection1.1$).
\end{prop}

\subsection{Locally algebraic abelian representations}
This subsection describes the essence of \cite[Ch. II, III]{Se98}.
Let $I_K$ be the id\`ele group of $K$, $\mathrm{Art}: I_K\to \Gal(\overline K/K)^{\mathrm{ab}}$ the \emph{Artin map}, and 
$i: (K\otimes_\Q\Q_\ell)^*\hookrightarrow I_K$ the embedding. Let $E_\ell$ be a finite extension of $\Q_\ell$ and
 $\bT:=\mathrm{Res}_{K/\Q}\mathbb{G}_m$ the $\Q$-torus given by the Weil restriction of scalars. 
For any $\Q$-algebra $A$, we have the natural group isomorphism $\bT(A)\cong (K\otimes_\Q A)^*$.

\begin{definition}\cite[Ch. III $\mathsection2$]{Se98}\label{loc}
A continuous abelian $E_\ell$-representation $\rho_\ell:\Gal(\overline K/K)^{\mathrm{ab}}\to \GL_n(E_\ell)$
is said to be \emph{locally algebraic}  if there exists an $E_\ell$-morphism 
$r:\bT\times_\Q E_\ell \to \GL_{n,E_\ell}$ such that
\begin{equation}
\rho_\ell\circ\mathrm{Art}\circ i(x)=r(x^{-1})
\end{equation}
for all $x\in(K\otimes_\Q \Q_\ell)^*$ close enough to the identity under the $\ell$-adic topology, 
where we identify $(K\otimes_\Q \Q_\ell)^*\subset (K\otimes_\Q E_\ell)^*=\bT(E_\ell)$.
\end{definition}

Local algebraicity of $\rho_\ell$ depends only on the local representations of $\rho_\ell$ above $\ell$, i.e.,
the restrictions of $\rho_\ell$ to the decomposition groups of $\Gal(\overline K/K)$ at the places of $\overline K$ above $\ell$.
This property is invariant under field extensions of $E_\ell$ (obviously) as well as the restriction of scalars.

\begin{prop}\label{loc-res}(see \cite[$\mathsection5$]{Ri76})
A continuous abelian $E_\ell$-representation $\rho_\ell:\Gal(\overline K/K)^{\mathrm{ab}}\to \GL_n(E_\ell)$ is 
locally algebraic if and only if it is locally algebraic when viewed as a $\Q_\ell$-representation.
\end{prop}

Let $S$ be a finite subset of finite places of $K$. For any \emph{modulus} $\mathfrak m=(m_v)_{v\in S}$ ($m_v\in\N$) 
with support $\mathrm{Supp}(\mathfrak{m})=S$,
 Serre has constructed \cite[Ch. II]{Se98} the \emph{Serre group} $\bS_\mathfrak{m}$ and an
 abelian rational representation with values in $\bS_{\mathfrak m}$ for each $\ell$:
\begin{equation}
\epsilon_\ell: \Gal(\overline K/K)^{\mathrm{ab}}\to\bS_\mathfrak{m}(\Q_\ell).
\end{equation}
The Serre group $\bS_\mathfrak{m}$ is a diagonalizable group over $\Q$ satisfying 
the short exact sequence
\begin{equation}\label{ses1}
1\to \bT_\mathfrak{m}\to \bS_\mathfrak{m}\to C_\mathfrak{m}\to 1,
\end{equation}
where $\bT_{\mathfrak{m}}$ is the identity component of $\bS_{\mathfrak{m}}$.

\begin{prop}\label{prop}\cite[Ch. II $\mathsection2$]{Se98} We list some properties of $\epsilon_\ell$ and $\bS_\mathfrak{m}$.
\begin{enumerate}[(i)]
\item The representation $\epsilon_\ell$ is unramified outside $S$ and $S_\ell$, the places $v$ of $K$ above $\ell$.
\item The image of the Frobenii above $v$ belongs to $\bS_\mathfrak{m}(\Q)$ for all $v\notin S\cup S_\ell$ 
and is dense in $\bS_\mathfrak{m}$.
\item The representations of $\bS_\mathfrak{m}$ are always semisimple. 
\item The finite group $C_{\mathfrak{m}}$ is the Galois group of the class field of $K$ corresponding to $\mathfrak{m}$. 
\item The dimension $d_K$ of $\bS_\mathfrak{m}$ depends only on $K$, not on $\mathfrak{m}$.
\end{enumerate}
\end{prop}

\begin{definition}\cite[Ch. II $\mathsection2$]{Se98}\label{Serre}
A continuous abelian $E_\ell$-representation $\rho_\ell:\Gal(\overline K/K)^{\mathrm{ab}}\to \GL_n(E_\ell)$
is said to be \emph{associated to some Serre group $\bS_\mathfrak{m}$} if there exists an $E_\ell$-morphism 
$\eta_\ell: \bS_{\mathfrak{m}, E_\ell} \to \GL_{n,E_\ell}$ such that 
\begin{equation}
\rho_\ell= \eta_\ell\circ \epsilon_\ell.
\end{equation}
\end{definition}

\begin{thm}\label{loc=Serre}\cite[Ch. III Thm. 1, Thm. 2]{Se98}
Let  $\rho_\ell:\Gal(\overline K/K)^{\mathrm{ab}}\to \GL_n(E_\ell)$ be a  continuous abelian $E_\ell$-representation.
\begin{enumerate}[(i)]
\item $\rho_\ell$ is locally algebraic if and only if it is associated to some Serre group $\bS_\mathfrak{m}$.
\item If $\rho_\ell$ is locally algebraic and unramified outside a finite set of places $S$, then 
it is associated to some Serre group $\bS_\mathfrak{m}$ with $\mathrm{Supp}(\mathfrak{m})\subset S$.
\end{enumerate}
\end{thm}

\begin{thm}\label{sc=>Sgp}\cite[Ch. I $\mathsection2.3$ Thm., Ch. II Thm. 1, Ch. III Thm. 2]{Se98}
Let $\{\rho_\ell\}_\ell$ be a rational strictly compatible system (Definition \ref{compat}) of abelian semisimple $\ell$-adic representations.
Then the system is associated to some Serre group $\bS_\mathfrak{m}$, i.e., there exists a $\Q$-morphism $\eta:\bS_\mathfrak{m}\to \GL_{n,\Q}$
such that $\rho_\ell= (\eta\otimes_\Q\Q_\ell)\circ \epsilon_\ell$ for all $\ell$.
\end{thm}

\begin{thm}\label{Serre=>rat}
If a continuous abelian $E_\ell$-representation $\rho_\ell:\Gal(\overline K/K)^{\mathrm{ab}}\to \GL_n(E_\ell)$
is associated to some Serre group $\bS_\mathfrak{m}$, then $\rho_\ell$ is semisimple and 
$E$-rational (Definition \ref{rat}) for some number field $E$.
\end{thm}

\begin{proof}
The semisimplicity follows from Proposition \ref{prop}(iii).
Let $S$ be the support of $\mathfrak{m}$, $F$ be the Galois closure of $K$, 
and the $|C_\mathfrak{m}|^{th}$ roots of unity in $\overline\Q$,
 where the finite group
$C_\mathfrak{m}$ is defined in (\ref{ses1}).
By \cite[Ch. II $\mathsection3.4$ Prop. 1]{Se98}, the roots of the characteristic polynomial $P_{v,\rho_\ell}(t)$ belong to $F$ for all $v\notin S\cup S_\ell$. This implies $P_{v,\rho_\ell}(t)\in (F\cap E_\ell)[t]$ for all $v\notin S\cup S_\ell$. 
Therefore, $\rho_\ell$ is $E$-rational for some number field $E\subset E_\ell$.
\end{proof}

\begin{cor}\label{Serre=>sys}
If a continuous abelian $E_{\ell'}$-representation $\rho_{\ell'}:\Gal(\overline K/K)^{\mathrm{ab}}\to \GL_n(E_{\ell'})$
is associated to some Serre group $\bS_\mathfrak{m}$, then there exist a number field $E$ and a strictly compatible $E$-rational system 
$\{\rho_\lambda\}_\lambda$ of semisimple abelian $\lambda$-adic representations with 
exceptional set contained in $\mathrm{Supp}(\mathfrak{m})$ and a local field isometry 
$E_{\lambda'}\cong E_{\ell'}$
for some $\lambda'$ such that 
\begin{equation}
\rho_{\ell'}=\rho_{\lambda'}.
\end{equation}
\end{cor}

\begin{proof}
By Theorem \ref{Serre=>rat} and \cite[Ch. II $\mathsection2.4$]{Se98}, $\eta_{\ell'}$ (in Definition \ref{Serre}) is the base change to $E_{\ell'}$ of an $E$-morphism $\eta:\bS_{\mathfrak{m}, E} \to \GL_{n,E}$ and we may assume the number field $E$ is dense in $E_{\ell'}$.
Then $\eta$ induces the desired system of $\lambda$-adic representations (see e.g. \cite[Ch. II Thm. 1]{Se98}).
Pick $\lambda'$ such that $E_{\lambda'}$ is isometric to $E_{\ell'}$, we obtain $\rho_{\ell'}=\rho_{\lambda'}$.
\end{proof}

The following \emph{converse} of Theorem \ref{Serre=>rat} (in light of Theorem \ref{loc=Serre}) follows from the arguments of Serre in \cite[Ch. III $\mathsection3$]{Se98} and
a result of transcendental numbers by Waldschmidt \cite{Wa81}.

\begin{thm}\label{rat=>loc}(see e.g. \cite[Thm. 2]{He80})
If a continuous semisimple abelian $E_\ell$-representation $\rho_\ell:\Gal(\overline K/K)^{\mathrm{ab}}\to \GL_n(E_\ell)$
is $E$-rational for some number field $E$, then $\rho_\ell$ is locally algebraic.
\end{thm}

\begin{definition}\label{HT}\cite[$\mathsection2$]{Se67}
Let $L$ be a finite extension of $\Q_\ell$ and $B_{HT}:=\bigoplus_{i\in\Z}\C_\ell(i)$ be the period ring, 
where the $\Gal(\overline L/L)$-module $\C_\ell(i)$ is the twist of $\C_\ell$ 
(the completion of $\overline L$) by the $i$th power of the cyclotomic character. 
A continuous $E_\ell$-representation $\varphi_\ell:\Gal(\overline L/L)\to \GL(V_\ell)$ 
is said to be \emph{of Hodge-Tate type} if the natural $\Gal(\overline L/L)$-inclusion
\begin{equation}\label{admissible}
B_{HT}\otimes_{L}(B_{HT}\otimes_{\Q_\ell}V_\ell)^{\Gal(\overline L/L)} \to B_{HT}\otimes_{\Q_\ell}V_\ell
\end{equation}
is an isomorphism.
\end{definition}

\begin{remark}
The $E_\ell$-representation $V_\ell$ is viewed as a $\Q_\ell$-representation in (\ref{admissible}).
\end{remark}

\begin{thm}\label{loc=HT}\cite{Ta67},\cite[Ch. III Appendix]{Se98}
A continuous abelian $\Q_\ell$-representation $\rho_\ell:\Gal(\overline K/K)^{\mathrm{ab}}\to \GL_n(\Q_\ell)$
is locally algebraic if and only if the local representations of 
$\varphi_\ell:\Gal(\overline K/K) \twoheadrightarrow\Gal(\overline K/K)^{\mathrm{ab}}\stackrel{\rho_\ell}{\rightarrow}\GL_n(\Q_\ell)$ above $\ell$ are 
\begin{enumerate}[(i)]
\item of Hodge-Tate type and
\item semisimple when restricted to the inertia groups.
\end{enumerate}
\end{thm}

Given a continuous semisimple abelian $E_\ell$-representation
\begin{equation}
\varphi_\ell:\Gal(\overline K/K)\twoheadrightarrow\Gal(\overline K/K)^{\mathrm{ab}}\stackrel{\rho_\ell}{\rightarrow}\GL_n(E_\ell),
\end{equation}
the results of this subsection are summarized by the diagram:
\begin{align}\label{equiv}
\begin{split}
\large
\xymatrix{
~\mathrm{B}~ \ar@{<=>}[d]_{\normalsize{\ref{loc=HT}}} \ar@{<=>}[r]^{\normalsize{\ref{loc-res}}} 
& ~\mathrm{A}~ \ar@{<=>}[d]_{\normalsize{\ref{loc=Serre}(i)}}  & ~\mathrm{D} ~ 
\ar@{=>}[l]_{\normalsize{\ref{rat=>loc}}} \ar@{<=}[d]^{\mathrm{obvious}} \ar@{<=}[ld]^{\normalsize{\ref{Serre=>rat}}}\\
~\mathrm{F}       ~               & ~\mathrm{C}~   \ar@{=>}[r]_{\normalsize{\ref{Serre=>sys}}}       & ~\mathrm{E}~,} 
\end{split}
\end{align}

where A-F are the following conditions and the arrows in \eqref{equiv} imply that they are equivalent.
\begin{enumerate}[A.]
\item $\rho_\ell$ is locally algebraic;
\item $\rho_\ell$ is locally algebraic as a $\Q_\ell$-representation;
\item $\rho_\ell$ is associated to some Serre group $\bS_\mathfrak{m}$;
\item $\rho_\ell$ is $E$-rational;
\item $\rho_\ell$ can be extended to a strictly compatible $E$-rational system of semisimple abelian $\lambda$-adic representations
with exceptional set contained in $\mathrm{Supp}(\mathfrak{m})$ for some modulus $\mathfrak m$;
\item The local representations above $\ell$ of $\varphi_\ell$ are of Hodge-Tate type.
\end{enumerate}

\begin{remark}
In later text, what we need from the diagram \eqref{equiv} are the equivalence $(\mathrm{F})\Leftrightarrow(\mathrm{C})$ and 
the implication $(\mathrm{C})\Rightarrow(\mathrm{E})$.
\end{remark}

\subsection{Formal bi-character}
Let $\bG_\C$ be a connected reductive Lie subgroup of $\GL_n(\C)$ and $\bG_\C^{\der}$ be the derived group (semisimple part) of $\bG_\C$.
Let $\bT_\C$ be a maximal torus of $\bG_\C$ and $\bT_\C^{\ss}:=\bT_\C\cap\bG_\C^{\der}$ be the maximal torus of $\bG_\C^{\der}$.
Since all maximal tori of $\bG_\C$ are conjugate (in $\GL_n(\C)$), the following definition makes sense.

\begin{definition}\label{fbc}\cite[Defs. 2.2, 2.3]{Hu18}
\begin{enumerate}[(i)]
\item The \emph{formal character} of $\bG_\C\hto\GL_n(\C)$ is defined as
the isomorphism class of the chain $\bT_\C\subset\GL_n(\C)$ up to $\GL_n(\C)$-conjugation. 
\item The \emph{formal bi-character} of $\bG_\C\hto\GL_n(\C)$ is defined as
the isomorphism class of the chain $\bT_\C^{\ss}\subset \bT_\C\subset\GL_n(\C)$ up to $\GL_n(\C)$-conjugation. 
\end{enumerate}
\end{definition}

For all $\ell$, let $E_\ell$ be a finite extension of $\Q_\ell$ and fix an embedding $E_\ell\hto\C$.
Consider a weakly compatible system of semisimple $E_\ell$-representations 
\begin{equation}\label{psi}
\psi_\ell:\Gal(\overline K/K)\to \GL_n(E_\ell),\hspace{.1in}\forall \ell
\end{equation}
and denote by $\bG_\ell$ the algebraic monodromy group of $\psi_\ell$.

\begin{thm}\label{thm4}
Suppose $\{\psi_\ell\}_\ell$ satisfies the following conditions:
\begin{enumerate}[(i)]
\item $\bG_\ell$ is connected for all $\ell$;
\item for each $\ell$, the eigenvalues of $\psi_\ell(\Frob_{\overline v})$ satisfy the conditions in Theorem \ref{Frob} for almost all finite places $v$ of $K$;
\item the local representations above $\ell$ of $\psi_\ell$ are of Hodge-Tate type for all $\ell$.
\end{enumerate}
Then the formal bi-character of $\bG_\ell(\C)\hto \GL_n(\C)$ is independent of $\ell$.
\end{thm}

\begin{remark}
In case $\{\psi_\ell\}_\ell$ is only a rational strictly compatible system of $\Q_\ell$-representations without the three conditions, Theorem \ref{thm4} follows from \cite[Thm. 3.19]{Hu13}. The added conditions allow us to streamline the proof of \cite[Thm. 3.19]{Hu13}, in which we follow closely for Theorem \ref{thm4}.
\end{remark}

\begin{proof}
Identify $\bG_\ell$ as a connected reductive subgroup of $\GL_{n,E_\ell}$ by (i) and the semisimplicity. Then the 
derived group $\bG_\ell^{\der}$ is connected and semisimple. 
Let $\eta:\bS_\mathfrak{m}\to\GL_{m,\Q}$ be a faithful representation (for any modulus $\mathfrak{m}$) 
and form a strictly compatible system 
of semisimple abelian $E_\ell$-representations
$\{\rho_\ell:=\eta_{\ell}\circ\epsilon_\ell\}_\ell$ (Definition \ref{Serre}). Then $\{\rho_\ell\}_\ell$ satisfies (ii) by \cite[Ch. II $\mathsection3.4$ Prop. 2]{Se98} and (iii) by the diagram (\ref{equiv}), so does the weakly compatible system
$\{\psi_\ell\oplus\rho_\ell\}_\ell$.
Since $\psi_\ell\oplus\rho_\ell$ is semisimple, 
its algebraic monodromy group $\widehat{\bG}_\ell$ is reductive.  
Let $Z(\widehat{\bG}_\ell^\circ)$ be the center of the identity component $\widehat{\bG}_\ell^\circ$.
Since $\rho_\ell$ is abelian, 
the derived group $\widehat{\bG}_\ell^{\der}$
of $\widehat{\bG}_\ell$ is exactly equal to $\bG_\ell^{\der}$ (or $\bG_\ell^{\der}\times\{1\}\subset\GL_{n,E_\ell}\times\GL_{m,E_\ell}$). Since $\widehat{\bG}_\ell$ is reductive,
there exist $h,k\in\Z_{\geq0}$ and a subrepresentation $W_\ell$ of
$$(\psi_\ell\oplus\rho_\ell)^{\otimes h}\otimes (\psi_\ell\oplus\rho_\ell)^{*\otimes k}$$
 such that 
\begin{equation}
\widehat{\bG}_\ell^{\der}=\mathrm{ker}(\widehat{\bG}_\ell\to \GL_{W_\ell})
\end{equation}
by \cite[I. Prop. 3.1(a)]{DMOS82}.
Since the category of Hodge-Tate type representations are stable under duality, tensor product, direct sum, and passage to submodules 
and quotients, the local representations above $\ell$ of $W_\ell$ is of Hodge-Tate type by construction.
Together with the fact that $W_\ell$ is abelian semisimple, $W_\ell$ is 
associated to some Serre group by the diagram (\ref{equiv}), which implies the inequality
\begin{equation}\label{dim1}
\dim \widehat{\bG}_\ell/\widehat{\bG}_\ell^{\der}
\leq \dim\bS_{\mathfrak{m},E_\ell}\stackrel{\small{\ref{prop}(v)}}{=:} d_K.
\end{equation}
Since $\widehat{\bG}_\ell^{\der}$ is a connected semisimple normal subgroup of $\widehat{\bG}_\ell^\circ$,
 the intersection of $\widehat{\bG}_\ell^{\der}$ and $Z(\widehat{\bG}_\ell^\circ)$ is finite.
Hence, we obtain
\begin{equation}\label{dim1.1}
\dim Z(\widehat{\bG}_\ell^\circ)\leq\dim \widehat{\bG}_\ell^\circ/\widehat{\bG}_\ell^{\der}
=\dim \widehat{\bG}_\ell/\widehat{\bG}_\ell^{\der}
\stackrel{\eqref{dim1}}{\leq}  d_K.
\end{equation}
Since the projection of $\widehat{\bG}_\ell\subset \GL_{n,E_\ell}\times\GL_{m,E_\ell}$ to the second factor is $\bS_{\mathfrak{m},E_\ell}$
by Proposition \ref{prop}(ii), we conclude with (\ref{dim1.1}) that 
\begin{equation}\label{dim2}
\dim Z(\widehat{\bG}_\ell^\circ)= d_K.
\end{equation}
We will see that this equality is crucial to establishing \eqref{goal} below.

Since $\{\psi_\ell\oplus\rho_\ell\}_\ell$ satisfies (ii) (the Frobenius eigenvalues conditions), 
the restriction $\{(\psi_\ell\oplus\rho_\ell)|_{\Gal(\overline K/L)}\}_\ell$ satisfies (ii) and weak compatibility
for any finite extension $L$ of $K$.
By Corollary \ref{connect}, we can find some finite extension $L$ such that the algebraic monodromy group of 
$(\psi_\ell\oplus\rho_\ell)|_{\Gal(\overline K/L)}$ is the identity component $\widehat{\bG}_\ell^\circ$ for all $\ell$
and the rank of $\widehat{\bG}_\ell^\circ$ is independent of $\ell$.
By Theorem \ref{Frob} and Remark \ref{FCremark} on 
the system $\{(\psi_\ell\oplus\rho_\ell)|_{\Gal(\overline K/L)}\}_\ell$, for a pair of primes $\ell_1,\ell_2$,
there exists a Frobenius element $\Frob_{\overline w}$ for some finite place $w$ of $L$ such that 
the semisimple element
\begin{align}\label{Frobss}
\begin{split}
\psi_{\ell_i}(\Frob_{\overline w})_{\ss}\times \rho_{\ell_i}(\Frob_{\overline w})_{\ss} 
\in \widehat{\bG}_{\ell_i}^\circ\subset\GL_n\times\GL_m
\end{split}
\end{align}
generates (algebraically) the maximal Frobenius tori $\widehat{\bT}_{\overline w,\ell_i}$ of $\widehat{\bG}_{\ell_i}^\circ$
for $i=1,2$.
Since the characteristic polynomial of $\psi_{\ell_i}(\Frob_{\overline w})_{\ss}$ 
(resp. $\rho_{\ell_i}(\Frob_{\overline w})_{\ss}$)
is independent of $i=1,2$,
the semisimple elements \eqref{Frobss} for $i=1,2$ are conjugate in $\GL_n(\C)\times\GL_m(\C)$.
By construction, the subtori
\begin{equation}\label{taut1}
\widehat{\bT}_{\overline w,\ell_i}(\C)\subset \GL_n(\C)\times\GL_m(\C)
\end{equation}
for $i=1,2$ are conjugate in $\GL_n(\C)\times\GL_m(\C)$.
Let $\widehat{\bT}_{\ell,\C}$ be a maximal torus of $\widehat{\bG}_{\ell}^\circ(\C)$ for all $\ell$.
Since $\ell_1$ and $\ell_2$ are arbitrary and all maximal tori of a reductive group over $\C$ are conjugate,
the subtorus
\begin{equation}\label{taut1.5}
\widehat{\bT}_{\ell,\C}\subset \GL_n(\C)\times\GL_m(\C)
\end{equation}
is independent of $\ell$ under conjugation by $\GL_n(\C)\times\GL_m(\C)$.

If we let $\pi_1$ (resp. $\pi_2$) be the projection of $\GL_n(\C)\times\GL_m(\C)$ to the first (resp. second) factor, then 
the $\ell$-independence of (\ref{taut1.5}) implies that the chain
\begin{equation}\label{taut2}
\bT^{\ss}_{\ell,\C}:=\pi_1(\mathrm{ker}(\pi_2|_{\widehat{\bT}_{\ell,\C}})^\circ)\subset
\bT_{\ell,\C}:=\pi_1(\widehat{\bT}_{\ell,\C})\subset \GL_n(\C)
\end{equation}
is independent of $\ell$ under conjugation by $\GL_n(\C)$.
Since $\pi_1(\widehat{\bG}_\ell)=\bG_\ell$, $\pi_1(\widehat{\bT}_{\ell,\C})$ is a maximal torus of $\bG_\ell(\C)$.
Since we have $\pi_2(\widehat{\bG}_\ell)=\bS_{\mathfrak{m},E_\ell}$ and (\ref{dim2}), we deduce that
\begin{equation}\label{goal}
\ker(\pi_2:\widehat{\bG}_\ell\to\GL_{m,E_\ell})^\circ=\widehat{\bG}_\ell^{\der}=\bG_\ell^{\der}\times\{1\}.
\end{equation}
Hence, $\pi_1(\mathrm{ker}(\pi_2|_{\widehat{\bT}_{\ell,\C}})^\circ)$ is a maximal torus of $\bG_\ell^{\der}(\C)$.
Therefore, we conclude that (\ref{taut2}), the formal bi-character of $\bG_\ell(\C)\hto\GL_n(\C)$, 
is independent of $\ell$.
\end{proof}

Suppose $\bG_\C$ and $\bG_\C'$ are two connected reductive Lie 
subgroups of $\GL_n(\C)$ such that their formal bi-characters are isomorphic. 
After conjugating $\bG_\C'$ by an element of $\GL_n(\C)$, we may assume 
\begin{equation}\label{cfbc}
\bT_\C^{\ss}\subset\bT_\C\subset\GL_n(\C)
\end{equation}
is the common formal bi-character of $\bG_\C$ and $\bG_\C'$,
meaning that $\bT_\C$ (resp. $\bT_\C^{\ss}$) 
is a common maximal torus of $\bG_\C$ and $\bG_\C'$ (resp. $\bG_\C^{\der}$ and $\bG_\C'^{\der}$). 
Let $\X^{\ss}$ be the character group of $\bT_\C^{\ss}$ and $R\subset\X^{\ss}$ (resp. $R'\subset\X^{\ss}$)
be the set of roots of $\bG_\C^{\der}$ (resp. $\bG_\C'^{\der}$) with respect to $\bT_\C^{\ss}$.

\begin{prop} \label{conjugation}\cite[Cor. 3.8]{Hu18}
Let $\bG_\C$ and $\bG_\C'$ be as above.
If $R=R'$ in $\X^{\ss}$, then $\bG_\C$ and $\bG_\C'$ are conjugate in $\GL_n(\C)$.
\end{prop}

We give sufficient conditions for $R=R'$ in $\X^{\ss}$.

\begin{prop}\label{R=R'}
Let $\bG_\C$ and $\bG_\C'$ be as above. Suppose the following conditions hold:
\begin{enumerate}[(i)]
\item $\bG_\C^{\der}$ and $\bG_\C'^{\der}$ are conjugate in $\GL_n(\C)$;
\item the set of roots $R$ of $\bG_\C^{\der}$ is stable under the normalizer $N_{\GL_n(\C)}(\bT_\C^{\ss})$.
\end{enumerate}
Then $R=R'$ in $\X^{\ss}$.
\end{prop}

\begin{proof} Since $\bG_\C^{\der}$ and $\bG_\C'^{\der}$ are conjugate in $\GL_n(\C)$ by (i), 
they are conjugate by an element of $N_{\GL_n(\C)}(\bT_\C^{\ss})$, inducing an 
automorphism of $\X^{\ss}$ mapping $R$ to $R'$. 
Hence, $R=R'$ by (ii).
\end{proof}

We return to the semisimple weakly compatible system (\ref{psi}) $\{\psi_\ell:\Gal(\overline K/K)\to\GL_n(E_\ell)\}_\ell$. 
The torus $\bT_{\ell,\C}^{\ss}$ in \eqref{taut2} is a maximal torus of $\bG_\ell^{\der}(\C)$. Let $R_\ell$ be
the set of roots of $\bG_\ell^{\der}(\C)$ with respect to $\bT_{\ell,\C}^{\ss}$.

\begin{cor}\label{l-indep}
Suppose $\{\psi_\ell\}_\ell$ satisfies the following conditions:
\begin{enumerate}[(i)]
\item $\bG_\ell$ is connected for all $\ell$;
\item for each $\ell$, the eigenvalues of $\psi_\ell(\Frob_{\overline v})$ satisfy the conditions in Theorem \ref{Frob} for almost all finite places $v$ of $K$;
\item the local representations above $\ell$ of $\psi_\ell$ are of Hodge-Tate type for all $\ell$;
\item the conjugacy class of $\bG_\ell^{\der}(\C)$ in $\GL_n(\C)$ is independent of $\ell$;
\item the set of roots $R_\ell$ of $\bG_\ell^{\der}(\C)$ is stable under the normalizer $N_{\GL_n(\C)}(\bT_{\ell,\C}^{\ss})$.
\end{enumerate}
Then the conjugacy class of $\bG_\ell(\C)$ in $\GL_n(\C)$ is independent of $\ell$.
In particular, the dimension of the abelian part $(E_\ell^n)^{\ab}$ (resp. the invariants $(E_\ell^n)^{\Gal(\overline K/K)}$)
is independent of $\ell$.
\end{cor}

\begin{proof}
The formal bi-character of $\bG_\ell(\C)\hto\GL_n(\C)$ is independent of $\ell$ by (i),(ii),(iii) and Theorem \ref{thm4}.
Then (iv),(v) and Propositions \ref{R=R'}, \ref{conjugation} imply that 
the conjugacy class of $\bG_\ell(\C)$ in $\GL_n(\C)$ is independent of $\ell$. 
\end{proof}

\subsection{Proofs of the theorems in $\mathsection1$}
$ $\\

\noindent \textit{\textbf{Proof of Theorem \ref{thm1}.}}
Since $\{V_\ell\}_\ell$ is a semisimple rational strictly compatible system arising from geometry, $\{V_\ell\}_\ell$
satisfies \ref{l-indep}(ii) and \ref{l-indep}(iii) respectively by Proposition \ref{eigenvalue} and 
the Hodge-Tate conjecture proved by Faltings \cite{Fa88}. 
Our conditions (i),(ii),(iii) are just \ref{l-indep}(i),(iv),(v). Hence, 
$\dim V_\ell^{\ab}$ is independent of $\ell$ by Corollary \ref{l-indep}.
Let $V_\ell^{\mathrm{c}}$ be the complementary representation of $V_\ell^{\ab}$ in $V_\ell$. 
Let $\chi_\ell$ be a one dimensional $\overline{\Q}_\ell$-representation (character) of $\Gal(\overline K/K)$. 
By the construction of $V_\ell^{\mathrm{c}}$, the twisted $\overline{\Q}_\ell$-representation
\begin{equation}\label{twist}
V_\ell^{\mathrm{c}}\otimes_{\Q_\ell}\chi_\ell~\mathrm{cannot~have~a~one~dimensional~subrepresentation.}
\end{equation} 

Pick a prime $\ell'$ and a finite extension $E_{\ell'}$ of $\Q_{\ell'}$ such that 
 $V_{\ell'}\otimes_{\Q_{\ell'}} E_{\ell'}$ decomposes as a direct sum of absolutely irreducible representations.
Consider the isomorphism classes $\mathcal{X}_{\ell'}$ of the characters appearing in $V_{\ell'}^{\ab}\otimes_{\Q_{\ell'}} E_{\ell'}$.  
For any $W_{\ell'}\in \mathcal{X}_{\ell'}$, $W_{\ell'}$ is unramified outside $S\cup S_{\ell'}$ (Definition \ref{compat}(i)).
Since $W_{\ell'}$ arises from geometry, its local representations above $\ell'$ are of Hodge-Tate type \cite{Fa88}.
By the diagram (\ref{equiv}) and Theorem \ref{loc=Serre}(ii), $W_{\ell'}$ is associated to some Serre group
$\bS_\mathfrak{m}$ with $\mathrm{Supp}(\mathfrak{m})\subset S\cup S_{\ell'}$. 
Since $\mathcal{X}_{\ell'}$ is finite, there exists a number field $E$ dense in $E_{\ell'}$ such that for every
$W_{\ell'}\in \mathcal{X}_{\ell'}$, $W_{\ell'}$ extends to a strictly compatible $E$-rational system of semisimple abelian 
$\lambda$-adic representations
$\{W_{\lambda}\}_{\lambda}$ with exceptional set contained in $S\cup S_{\ell'}$ by the diagram (\ref{equiv}).
Suppose $\{W_{1,\lambda}\}_\lambda,\{W_{2,\lambda}\}_\lambda$ are the $\lambda$-adic systems 
corresponding to $W_{1,\ell'},W_{2,\ell'}\in\mathcal{X}_{\ell'}$. 
By the strict compatibility of the systems, 
\begin{equation}\label{separate}
\mathrm{if}~ W_{1,\ell'}\ncong W_{2,\ell'},~\mathrm{ then}~ W_{1,\lambda}\ncong W_{2,\lambda}~\mathrm{for~all}~\lambda.
\end{equation}
Let $\lambda'$ be the prime of $E$ corresponding to $E\hto E_{\ell'}$. 
For every prime $\ell$ other than $\ell'$, fix a $\lambda$ that divides $\ell$ and let $E_\ell$ be the completion of $E$ with 
respect to $\lambda$. Then we obtain for every 
$W_{\ell'}\in \mathcal{X}_{\ell'}$, a subfamily $\{W_{\ell}\}_\ell$ of $\{W_{\lambda}\}_{\lambda}$, indexed by $\ell$,
satisfying the following properties.\\

\noindent ($*$): \textit{The system $\{W_{\ell}\}_\ell$ extends $W_{\ell'}$ and is a weakly compatible system of $E_\ell$-representations. Moreover,  $W_\ell$ is unramified outside $(S\cup S_{\ell'})\cup S_\ell$.}\\

\noindent The invariant dimension of the representation 
$\mathrm{Hom}_{E_\ell}(W_\ell,V_\ell\otimes_{\Q_\ell} E_\ell)=V_\ell\otimes_{\Q_\ell} W_\ell^*$, that is,
\begin{equation}\label{inv}
\dim_{E_\ell}(V_\ell\otimes_{\Q_\ell} W_\ell^*)^{\Gal(\overline K/K)},
\end{equation}
counts the number of copies of $W_\ell$ in $V_\ell^{\ab}\otimes E_\ell$ by (\ref{twist}).
Temporarily assuming that (\ref{inv}) is independent of $\ell$ for each $W_{\ell'}\in \mathcal{X}_{\ell'}$,
then the facts that $\{W_\ell\}_\ell$ extends $W_{\ell'}$ and $\dim V_\ell^{\ab}$ is independent of $\ell$, together with
(\ref{separate}) and ($*$),
imply that for any pair of primes $\ell'',\ell'''$, 
\begin{equation}\label{charpoly}
P_{v,V_{\ell''}^{\ab}}(t)=P_{v,V_{\ell'''}^{\ab}}(t)\in E[t]~\mathrm{if~}v\notin (S\cup S_{\ell'})\cup S_{\ell''}\cup S_{\ell'''}.
\end{equation}
Roughly speaking, it says that the eigenvalues of Frobenius on $V_\ell^{\ab}$ is independent of $\ell$.
Let $P_v(t)$ be the characteristic polynomial in (\ref{charpoly}). Then $P_v(t) \in\Q_\ell[t]$ for almost all $\ell$ because $V_\ell^{\ab}$ is a $\Q_\ell$-vector space. Since $P_v(t)\in E[t]$ and $E$ is a finite extension of $\Q$, we obtain $P_v(t)\in\Q[t]$. 
Thus, $\{V_\ell^{\ab}\}_\ell$ is a rational strictly compatible with exceptional set contained in  $S\cup S_{\ell'}$.
By picking another prime $\ell''$ not equal to $\ell'$ in the second paragraph of this proof and repeat the above arguments, we conclude 
that $\{V_\ell^{\ab}\}_\ell$ is a rational strictly compatible with exceptional set contained in $S$.

To complete the proof, it remains to prove that (\ref{inv}) is independent of $\ell$ for each $W_{\ell'}\in \mathcal{X}_{\ell'}$.
It suffices to verify the five conditions of Corollary \ref{l-indep} for the semisimple weakly compatible system 
$\{V_\ell\otimes_{\Q_\ell} W_\ell^*\}_\ell$ for any $W_{\ell'}\in\mathcal{X}_{\ell'}$. 
Let $\bH_\ell$ be the algebraic monodromy group of $V_\ell\otimes_{\Q_\ell} W_\ell^*=(V_\ell\otimes_{\Q_\ell} E_\ell) \otimes_{E_\ell} W_\ell^*$. Since $\bG_\ell$ acts naturally on $W_\ell^*$,
the natural morphism 
$$\bG_\ell\to \GL_{(V_\ell\otimes_{\Q_\ell} E_\ell) \otimes_{E_\ell} W_\ell^*}$$
 implies that $\bG_\ell$ surjects 
onto $\bH_\ell$. Since $\bG_\ell$ is connected by (i), $\bH_\ell$ is also connected and \ref{l-indep}(i) holds.
Since $V_\ell\otimes_{\Q_\ell} W_\ell^*$ and $V_\ell\otimes_{\Q_\ell}E_\ell$
differ by a twist of the character $W_\ell^*$, 
the derived subgroups of the Galois images are equal under suitable identification, which implies
$\bG_\ell^{\der}(\C)=\bH_\ell^{\der}(\C)$ in $\GL_N(\C)$.
The conditions \ref{l-indep}(iv),(v) then follow from our (ii),(iii).
Since $\{W_\ell\}_\ell$ arises from an $E$-representation of some Serre group $\bS_{\mathfrak{m},E}$, 
$\{W_\ell\}_\ell$ satisfies \ref{l-indep}(ii) by \cite[Ch. II $\mathsection3.4$]{Se98} and \ref{l-indep}(iii) by the diagram (\ref{equiv}).
Hence, $V_\ell\otimes_{\Q_\ell} W_\ell^*$ also satisfies \ref{l-indep}(ii),(iii).
We conclude by Corollary \ref{l-indep} that the invariant dimension (\ref{inv}) is independent of $\ell$ for each $W_{\ell'}\in\mathcal{X}_{\ell'}$. \qed

\begin{thm}\label{hypoA}\cite[$\mathsection3$]{Hu18}
Let $\bG_\C$ and $\bG_\C'$ be two connected reductive Lie subgroup of $\GL_N(\C)$ having the 
common formal bi-character $\bT_\C^{\ss}\subset\bT_\C\subset\GL_N(\C)$. 
Suppose  $\bG_\C'^{\der}$ satisfies the Lie type conditions in Hypothesis A, i.e.,
\begin{itemize}
\item $\bG_\C'^{\der}$ has at most one factor of type $A_4$;
\item if $\bH_\C$ is an almost simple factor of $\bG_\C'^{\der}$, then 
$\bH_\C$ is of type $A_n$ for some
$n\in\N\backslash\{1,2,3,5,7,8\}$.
\end{itemize}
Then $\bG_\C$ and $\bG_\C'$ are conjugate in $\GL_N(\C)$ and their roots with respect to 
$\bT_\C$ are identical. In particular,
$\bG_\C^{\der}$ and $\bG_\C'^{\der}$ are conjugate in $\GL_N(\C)$ and their roots with respect to 
$\bT_\C^{\ss}$ are identical.
\end{thm}

\noindent \textit{\textbf{Proof of Theorem \ref{thm2}.}}
(i+ii): It suffices to show that Hypothesis A implies Theorem \ref{thm1}(ii),(iii).
This is true by Theorem \ref{thm4} and Theorem \ref{hypoA}.

(i+ii'): It follows directly from Theorems \ref{thm1} and \ref{thm6}.\qed\\

\noindent \textit{\textbf{Proof of Theorem \ref{thm3}.}}
The idea is similar to Theorem \ref{thm1}. Consider
the semisimple weakly compatible system 
\begin{equation}\label{newsys}
\{\mathrm{Hom}_{E_\ell}(W_\ell,V_\ell)=V_\ell\otimes_{\Q_\ell}W_\ell^*\}_\ell.
\end{equation}
Let $\bH_\ell$ be the algebraic monodromy group of  $V_\ell\otimes_{\Q_\ell}W_\ell^*$.
As $\{W_\ell\}_\ell$ is absolutely irreducible extending $W_{\ell'}$, 
it suffices to show that invariant dimension 
\begin{equation}
\dim_{E_\ell}(V_\ell\otimes_{\Q_\ell}W_\ell^*)^{\Gal(\overline K/K)}
\end{equation}
is independent of $\ell$, which follows from Corollary \ref{l-indep} once we 
check the five conditions for the system (\ref{newsys}).

The natural morphism 
$$\bG_\ell\to \GL_{(V_\ell\otimes_{\Q_\ell}E_\ell)\otimes_{E_\ell}W_\ell^*}=\GL_{V_\ell\otimes_{\Q_\ell}W_\ell^*}$$
 induces the 
surjection $\bG_\ell\to\bH_\ell$. Since $\bG_\ell$ is connected for all $\ell$ (i), $\bH_\ell$ is also connected 
and we obtain \ref{l-indep}(i).
Since $W_{\ell'}$ is arises from geometry and $\{W_\ell\}_\ell$ is weakly compatible satisfying (iv), the system 
(\ref{newsys}) satisfies \ref{l-indep}(ii),(iii). 
Theorem \ref{thm4} then implies that 
\begin{equation}\label{Afb}
\text{the formal bi-character of~} \bH_\ell(\C)\hto\GL_M(\C)\mathrm{~ is~ independent~ of~} \ell
\end{equation}
if $M$ is the dimension of (\ref{newsys}). 
Since $\{\bG_\ell\}_\ell$ satisfies Hypothesis A and $\bG_\ell$ surjects onto $\bH_\ell$ for all $\ell$,
the family $\{\bH_\ell(\C)\}_\ell$ of connected reductive Lie groups also satisfies Hypothesis A.
Therefore, \ref{l-indep}(iv),(v) follow from Theorem \ref{hypoA} and (\ref{Afb}).\qed

\subsection{$\ell$-independence of the Tate conjecture}
$ $\\

\noindent \textit{\textbf{Proof of Corollary \ref{indepTate}.}}
Let $X$ be a smooth projective variety defined over $K$ and $V_\ell:=H^{2r}(\overline{X},\Q_\ell(r))$.
Let $\mathfrak{g}_\ell$ be the Lie algebra of the image of $\Gal(\overline K/K)$ in $\GL(V_\ell)$ and $CH^r(\overline{X})$ be
the Chow group of codimension $r$ cycles of $\overline{X}$. Embed $\overline K$ in $\C$. 
The diagram
\begin{equation}\label{Chdiagram}
\xymatrix{
CH^r(\overline{X}) \ar[d] \ar[r]^{c_\ell} & H^{2r}(\overline{X},\Q_\ell(r))\ar[d]^{\alpha_\ell\cong}\\
CH^r(X(\C)) \ar[r]^{c_\C\otimes 1\hspace{.2in}} & H^{2r}(X(\C),\Q)\otimes\Q_\ell}
\end{equation}
is commutative, where $\alpha_\ell$ is the comparison isomorphism between \'etale and singular cohomology 
and $c_\C$ is the cycle map
for the complex manifold $X(\C)$. 
Suppose the algebraic monodromy group $\bG_\ell$ is connected for all $\ell$ by taking a large enough extension of $K$.
Then by semisimplicity we obtain
\begin{equation}\label{Tateinv}
V_\ell^{\mathfrak{g}_\ell}=V_\ell^{\bG_\ell}=V_\ell^{\Gal(\overline K/K)}=(V_\ell^{\mathrm{ab}})^{\Gal(\overline K/K)}.
\end{equation}
By Conjecture \ref{A} and Theorem \ref{sc=>Sgp}, the dimension of (\ref{Tateinv}) is independent of $\ell$.
Since $\mathrm{Im}c_\ell\subset V_\ell^{\mathfrak{g}_\ell}$ and 
the dimension of $\Q_\ell\cdot\mathrm{Im}c_\ell$ is 
independent of $\ell$ by (\ref{Chdiagram}), we conclude that the Tate conjecture is independent of $\ell$.\qed

\begin{remark}\label{bigenough}
The above proof only needs Conjecture \ref{A} for a big enough extension of $K$.
\end{remark}

\noindent \textit{\textbf{Proof of Corollary \ref{TateAV}.}}
Let $r\geq 0$ be an integer. Since the Galois representation $H^{1}(\overline{X},\Q_\ell)$ is semisimple for all $\ell$ \cite{Fa83},
so is the representation $V_\ell:=H^{2r}(\overline X,\Q_\ell(r))\cong\bigwedge^{2r}H^{1}(\overline{X},\Q_\ell)(r)$ for all $\ell$.
Let $\bH_\ell$ (resp. $\bG_\ell'$) be the algebraic monodromy group of $V_\ell$ 
(resp. $\bigwedge^{2r}H^{1}(\overline{X},\Q_\ell)$) for all $\ell$. 
By Corollary \ref{indepTate}, Theorem \ref{thm1}, and Remark \ref{bigenough}, 
it suffices to show that the groups $\bH_\ell$ satisfy the condition (ii) or (ii') of Theorem \ref{thm2}.
Since the groups $\bG_\ell$ satisfy \ref{thm2}(ii) or (ii'), the groups $\bG_\ell'$ also satisfy \ref{thm2}(ii) or (ii') by construction.
Since $\bigwedge^{2r}H^{1}(\overline{X},\Q_\ell)$ and $V_\ell$ are differed by a twist of character, the tautological representations 
$(\bG_\ell')^{\der}\hookrightarrow\GL_{\bigwedge^{2r}H^{1}}$ and $\bH_\ell^{\der}\hookrightarrow\GL_{V_\ell}$ are isomorphic for all $\ell$.
Therefore, the groups $\bH_\ell$ satisfy \ref{thm2}(ii) or (ii') and we are done.\qed

\appendix
\section{Geometry of roots}
Let $\bG$ be a connected semisimple complex Lie subgroup of $\GL_N(\C)$ and $\bT$ be a maximal torus of $\bG$.
Denote by $\X$ the character group of $\bT$ and by $R\subset\X$ the set of roots of $\bG$ with respect to $\bT$.
Let $f:\bG\hto\GL_N(\C)$ be the tautological representation. Then
the normalizer $N_{\GL_N(\C)}(\bT)$ acts on $\X$. We are interested in the following question:\\

\noindent\textbf{Q}: For what pair $(\bG,f)$ is $R$ stable under $N_{\GL_N(\C)}(\bT)$?\\

Theorem \ref{thm5} (a consequence of Theorem \ref{hypoA}) below answers \textbf{Q} for some type A semisimple $\bG$. 
The goal of this appendix is to study \textbf{Q} under the assumption that $\bG$ is almost simple (Theorem \ref{thm6}),
i.e., the Lie algebra $\mathfrak{g}$ of $\bG$ is a complex simple Lie algebra. In the Cartan-Killing classification,
$\mathfrak{g}$ is one of $A_n~(n\geq 1)$, $B_n~(n\geq 2)$, $C_n~(n\geq 3)$,
$D_n~(n\geq 4)$, $E_6$, $E_7$, $E_8$, $F_4$, and $G_2$.

\begin{thm}\label{thm5}
Let $\bG$ be a connected semisimple complex Lie group satisfying the Lie type conditions in Hypothesis A, i.e.,
\begin{itemize}
\item $\bG$ has at most one factor of type $A_4$;
\item if $\bH$ is an almost simple factor of $\bG$, then 
$\bH$ is of type $A_n$ for some
$n\in\N\backslash\{1,2,3,5,7,8\}$.
\end{itemize}
Then $R$ is stable under $N_{\GL_N(\C)}(\bT)$ for all embedding $f$ of $\bG$.
\end{thm}

\begin{thm}\label{thm6}
Let $\bG$ be a connected almost simple complex Lie group.
\begin{enumerate}[(i)]
\item If the type of $\bG$ is different from $A_7,A_8,B_4,D_8$, then $R$ is 
stable under $N_{\GL_N(\C)}(\bT)$ for all embeddings $f$ of $\bG$.
\item For all $\mathfrak{g}\in\{A_7,A_8,B_4,D_8\}$, there exists a pair $(\bG,f)$ such that $\bG$ 
is of type $\mathfrak{g}$ and $R$ is not stable under $N_{\GL_N(\C)}(\bT)$.
\end{enumerate}
\end{thm}

\subsection{A geometric viewpoint}
We introduce a geometric viewpoint 
on the configuration of $R$ in $\X\otimes_\Z\R$ that we learned from \cite[$\mathsection1$]{LP90}
and developed in \cite[$\mathsection2$]{Hu13}.
This viewpoint is also the starting point of Theorem \ref{hypoA}, see \cite[$\mathsection2$]{Hu13},\cite[$\mathsection3$]{Hu18}.
Let $\X$ be the character group of the maximal torus $\bT$, $\Z[\X]$ be the group ring, and $R$ be
(resp. $W$) the set of roots (resp. the Weyl group) of the connected complex semisimple $\bG$ with respect to $\bT$.
Then the formal character $\bT\subset\GL_N(\C)$ of 
$f:\bG\hto\GL_N(\C)$ (Definition \ref{fbc}) is equivalent to the element
\begin{equation}\label{Charf}
\Char(f):=\alpha_1+\cdots+\alpha_N\in\Z[\X],
\end{equation} 
where $\alpha_1,...,\alpha_N$ are the multiset of the weights of $f$ with respect to $\bT$.
Consider the $\R$-vector space $\X\otimes_\Z\R$.
The formal character defines a positive definite 
inner product $(~,~)$ on the dual space $(\X\otimes_\Z\R)^*$ by
\begin{equation}
(x_1^*,x_2^*)=\sum_{i=1}^N\alpha_i(x_1^*)\alpha_i(x_2^*). 
\end{equation}
This non-degenerate pairing induces an isomorphism $\X\otimes_\Z\R\to(\X\otimes_\Z\R)^*$ mapping $x$ to $x^*$.
Then we obtain a positive definite inner product $\left\langle~, ~\right\rangle$ on $\X\otimes_\Z\R$
such that $\left\langle x_1, x_2\right\rangle= (x_1^*,x_2^*)$ for all $x_1,x_2\in \X\otimes_\Z\R$.
Since the normalizer $\bN:=N_{\GL_N(\C)}(\bT)$ fixes $\Char(f)$,
$\bN$ is isometric on the metric space $\X\otimes_\Z\R$ induced by $\left\langle~, ~\right\rangle$.
Since the Weyl group $W$ is a quotient of $\bN_{\bG}(\bT)$, we conclude that
$W$ is isometric on $\X\otimes_\Z\R$.

Consider the set of almost simple factors of $\bG$:  $\bH_{j}$ for all $1\leq j\leq m$.
Then $\bT_{j}:=\bT\cap \bH_{j}$ is a maximal torus of $\bH_{j}$.
Let $\X_j$ be the character group of $\bT_{j}$ and $R_j$ (resp. $W_j$) be the set of roots (resp. Weyl group) 
of $\bH_{j}$ with respect to the maximal torus $\bT_{j}$.
Since $R_j$ spans $\X_j\otimes_\Z\R$,
the disjoint union
$R=\coprod_{j=1}^m R_j$
induces a direct sum decomposition:
\begin{equation}\label{rootsp}  
\X\otimes_\Z\R=\bigoplus_{j=1}^m\X_j\otimes_\Z\R.
\end{equation}
As $W\cong \prod_{j=1}^m W_j$ is isometric on $\X\otimes_\Z\R$ and preserves $\X_j\otimes_\Z\R$ for all $1\leq j\leq m$,
(\ref{rootsp}) is an orthogonal direct sum. This implies by \cite[VI $\mathsection1$ Prop. 5 Cor. (i)]{Bo81} that for $1\leq j\leq m$,
\begin{equation}\label{Killing}
\left\langle~, ~\right\rangle|_{\X_j\otimes_\Z\R}~ \text{is a positive multiple of the \textit{Killing form} \cite[$\mathsection14$]{FH91}
on} ~\X_j\otimes_\Z\R.
\end{equation}

Suppose $\bG'$ is another connected complex semisimple Lie subgroup of $\GL_N(\C)$ 
sharing the same maximal torus $\bT$ with $\bG$.
Let $R'$ be the set of roots of $\bG'$ with respect to $\bT$. Then $R'\subset\X\subset \X\otimes_\Z\R$.
 There is a strong geometric connection between $R$ and $R'$.
Let $r\in R$, $r'\in R'$, and $\theta$ be the angle between them.  
Since $r'$ and $r$ are characters of $\bT$, we obtain 
\begin{equation}\label{angle}
2\frac{\left\langle r', r\right\rangle}{\left\langle r, r\right\rangle} \in \Z, \hspace{.2in} 2\frac{\left\langle r, r'\right\rangle}{\left\langle r', r'\right\rangle} \in \Z
\end{equation}
by (\ref{Killing}). This implies $4\cos^2\theta\in\Z$ and $\theta$ can only be 
a multiple of $30^\circ$ or $45^\circ$ (see \cite[$\mathsection2$]{Hu13}).\\

\begin{minipage}{0.6\textwidth}
    \noindent \textbf{Example}: Notation as above. The semisimple Lie group
		$\bG=\SL_2(\C)\times\SL_2(\C)/\pm(-\mathrm{Id},-\mathrm{Id})$ is 
a subgroup of the simple exceptional Lie group $\bG'$ of type $G_2$ and $f$ is the restriction to $\bG$ of 
the adjoint representation $\mathrm{Ad}:\bG'\hto\GL_N(\C)$. Since $\bT$ is also a maximal torus of 
$\bG'$, the configuration of the nonzero weights of $f$ on $\X\otimes_\Z\R$ with respect to $\left\langle~, ~\right\rangle$
is just some positive multiple of the root system of $G_2$ (the figure).
As $f$ contains the adjoint representation of $\bG$ as a subrepresentation, the set of roots $R=R_1\coprod R_2$ of $\bG$ must appear 
in the figure. Since (\ref{rootsp}) is an orthogonal sum, $R_1$ and $R_2$ are perpendicular. 
Without loss of generality, assume $R_1=\{e_1,-e_1\}$ and $R_2=\{e_2,-e_2\}$. Then we obtain 
$||e_2||=\sqrt{3}||e_1||$.
On the other hand, we see that the angles between $R$ (red) and $R'$ (red or blue)
are multiples of $30^\circ$.
\end{minipage}	\hfill
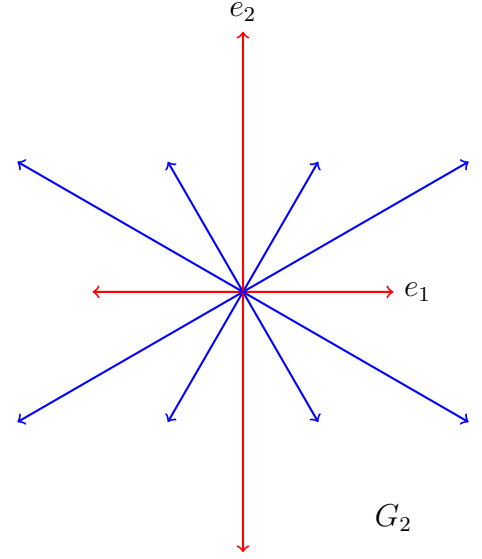
\begin{minipage}{.1\textwidth}
\end{minipage}	\hfill
\begin{minipage}{0.35\textwidth}  
\begin{tikzpicture}
		\draw[->,thick, red] (0, 0) -- (0:2);
		\draw[->,thick, blue] (0, 0) -- (60:2);
		\draw[->,thick, blue] (0, 0) -- (2*60:2);
		\draw[->,thick, red] (0, 0) -- (3*60:2);
		\draw[->,thick, blue] (0, 0) -- (4*60:2);
		\draw[->,thick, blue] (0, 0) -- (5*60:2);
		\draw[->,thick, blue] (0, 0) -- (30 + 0*60:3.464);
		\draw[->,thick, red] (0, 0) -- (30 + 1*60:3.464);
		\draw[->,thick, blue] (0, 0) -- (30 + 2*60:3.464);
		\draw[->,thick, blue] (0, 0) -- (30 + 3*60:3.464);
		\draw[->,thick, red] (0, 0) -- (30 + 4*60:3.464);
		\draw[->,thick, blue] (0, 0) -- (30 + 5*60:3.464);
    \node[right] at (2, 0) {$e_1$};
		\node[above] at (0, 3.464) {$e_2$};
    \node at (2, -3) {$G_2$};
  \end{tikzpicture}
\end{minipage}

\subsection{Proof of Theorem \ref{thm6}}
The notation in $\mathsection A.1$ remains in force. 
There is a lattice $\Lambda$ of $\X\otimes_\Z\R$ containing $\X$, known as the \emph{weight lattice}.
The construction is as follows. Let $\widetilde{\bG}$ be the universal covering of $\bG$, $\widetilde{\bT}$ be
a maximal torus of $\widetilde{\bG}$ mapping onto $\bT$, and $\Lambda$ be the character group of $\widetilde{\bT}$.
Then the character group $\X$ of $\bT$ injects into $\Lambda$ as a finite index subgroup. Hence, $\Lambda\otimes_\Z\R$
can be identified as  $\X\otimes_\Z\R$. 
Suppose $\bG$ is simple and let $g\in\bN$.
Then $g$ is orthogonal on $\Lambda\otimes_\Z\R$.
Denote by $\bG'$ the subgroup $g\bG g^{-1}$ in $\GL_N(\C)$.
Then $\bG$ and $\bG'$ share the common maximal torus $\bT$ and $R'=g R$.
We first prove that $R=R'$ in $\Lambda$ if the Lie algebra $\mathfrak{g}$ of $\bG$ is 
not equal to $A_7,A_8,B_4,D_8$ (Theorem \ref{thm6}(i)).
The configuration $R\subset\Lambda\subset\Lambda\otimes_\Z\R$, up to normalization, is the same as
the configuration of the root system of the complex simple Lie algebra $\mathfrak{g}$ \cite[Table 1]{GOV94}.

\subsubsection{$\mathfrak{g}=A_n$}
Let $r\in R$, $r'\in R'$, and $\theta$  be the angle between them.
Since $R'=g R$, we obtain $||r||=||r'||$. Thus, 
$\theta$ is a multiple of $60^\circ$ or $90^\circ$ ($2\cos\theta\in \Z$ by (\ref{angle})).
We follow the crucial computation in \cite[Prop. 2.8]{Hu13}.
Normalize $\Lambda\otimes_\Z\R$ such that $||r||=\sqrt{2}$. 
There exist weights $e_1,...,e_n,e_{n+1}\in \Lambda$ such that 
\begin{align}\label{A_n}
\begin{split}
0 &= e_1+e_2+...+e_n+e_{n+1}\\
\Lambda &=\mathbb{Z}e_1\oplus\cdots\oplus\mathbb{Z}e_n\\
\left\langle e_i, e_j \right\rangle &= \left\{ \begin{array}{lll}
 \frac{n}{n+1} &\mbox{if}& 1\leq i=j\leq n\\
 \frac{-1}{n+1} &\mbox{if}& 1\leq i\neq j\leq n
\end{array}\right.
\end{split}
\end{align}
and the set of roots $R$ comprise the set $\{e_i-e_j:$ $1\leq i\neq j\leq n+1\}$. 
Since $r'\in\Lambda$, write $r'=a_1e_1+\cdots+a_ne_n$ as an integral combination.
Consider $r=e_i-e_j$ where $i,j\leq n$. By (\ref{A_n}), we have 
\begin{align}\label{aij}
\begin{split}
2\cos\theta=\left\langle r', r \right\rangle &=\left\langle a_1e_1+\cdots +a_ne_n, e_i-e_j \right\rangle\\
&=\frac{a_i(n+1)-\sum_{k=1}^na_k}{n+1}-\frac{a_j(n+1)-\sum_{k=1}^na_k}{n+1}=a_i-a_j.
\end{split}
\end{align}
Then consider $r=e_i-e_{n+1}$ where $i\leq n$. By (\ref{A_n}), we have $r=e_1+e_2+\cdots +e_{i-1}+2e_i+e_{i+1}+\cdots +e_n$, so
\begin{align}\label{ai}
\begin{split}
2\cos\theta=\left\langle r', r \right\rangle &=\left\langle a_1e_1+\cdots+a_ne_n, e_1+e_2+\cdots+e_{i-1}+2e_i+e_{i+1}+\cdots+e_n \right\rangle\\
&=\frac{a_i(n+1)-\sum_{k=1}^na_k}{n+1}+\sum_{j=1}^n\frac{a_j(n+1)-\sum_{k=1}^na_k}{n+1}=a_i.
\end{split}
\end{align}
Suppose $r'\notin R$, then $\theta\in\{60^\circ,90^\circ,120^\circ,240^\circ,270^\circ,300^\circ\}$ and 
\begin{equation}
|a_i|,|a_i-a_j| \in\{0,1\},\hspace{.1in}\forall 1\leq i\neq j\leq n
\end{equation}
by (\ref{aij}) and (\ref{ai}), which means that either all $a_i\in\{0,1\}$ or all $a_i\in\{0,-1\}$. 
Therefore, we obtain
\begin{align}\label{Dio}
\begin{split}
2=\left\langle r', r' \right\rangle &= \left\langle a_1e_1+\cdots+a_ne_n, a_1e_1+\cdots+a_ne_n \right\rangle=\frac{\sum_{i=1}^nna_i^2-2\sum_{i<j}a_ia_j}{n+1}\\
&=\frac{\sum_{i=1}^na_i^2+\sum_{i< j}(a_i-a_j)^2}{n+1}=\frac{(n-k)+k(n-k)}{n+1},
\end{split}
\end{align}
where $k$ is the number of zero $a_i$.
Since the integral equation (\ref{Dio}) has a solution $(k,n)\in\Z_{\geq 0}\times\N$ 
such that $k\leq n$ (i.e., $r'\notin R$) 
if and only if $n=7,8$ \cite[(\textbf{2.13}) $\theta=60^\circ$]{Hu13},
we conclude that $R=R'$ if $n\in\N\backslash\{7,8\}$.

\subsubsection{$\mathfrak{g}=B_n$}
Under suitable normalization, $\Lambda\otimes_\Z\R$ is isometric to the Euclidean space $\R^n$. 
The weight lattice $\Lambda$ is the lattice generated by the standard basis $\{e_1,...,e_n\}$ and 
$\frac{e_1+\cdots +e_n}{2}$. The set of roots $R=\{\pm e_i, \pm e_i\pm e_j:~1\leq i\neq j\leq n\}$.

Let $r'\in R'$ be a short root, i.e., $||r'||=1$. Write $r'=\frac{a_1e_1+\cdots+a_ne_n}{2}$ as
a half-integral combination. Then either all $a_i$ are even or all $a_i$ are odd.
In the even case, $||r'||=1$ implies $r'\in R$. In the odd case, $||r'||=1$ implies $n=4$.
Hence, the short roots $R'^\circ$ of $R'$ belong to $R$ whenever $n\neq 4$. 
Since $R'^\circ$ determines $R'$, we obtain $R=R'$ 
if $n\neq 4$. 

\subsubsection{$\mathfrak{g}=C_n$}
Under suitable normalization, $\Lambda\otimes_\Z\R$ is isometric to the Euclidean space $\R^n$. 
The weight lattice $\Lambda$ is the lattice generated by the standard basis $e_1,...,e_n$.
The set of roots $R=\{\pm 2e_i, \pm e_i\pm e_j:~1\leq i\neq j\leq n\}$.

Since all weights of length $\sqrt{2}$ belong to $R$, 
the short roots $R'^\circ$ of $R'$ belong to $R$. Since $R'^\circ$ determines $R'$, we obtain $R=R'$.

\subsubsection{$\mathfrak{g}=D_n$}
Under suitable normalization, $\Lambda\otimes_\Z\R$ is isometric to the Euclidean space $\R^n$. 
The weight lattice $\Lambda$ is the lattice generated by the standard basis $e_1,...,e_n$ and 
$\frac{e_1+\cdots +e_n}{2}$ (same as $B_n$).
The set of roots $R=\{\pm e_i\pm e_j:~1\leq i\neq j\leq n\}$.

Let $r'\in R'$ be a root, $||r'||=\sqrt{2}$. Write $r'=\frac{a_1e_1+\cdots+a_ne_n}{2}$ as
a half-integral combination. Then either all $a_i$ are even or all $a_i$ are odd.
In the even case, $||r'||=\sqrt{2}$ implies $r'\in R$. In the odd case, $||r'||=\sqrt{2}$ implies $n=8$.
This implies  $R=R'$ if $n\neq 8$. 

\subsubsection{$\mathfrak{g}=E_6,E_7,E_8$} Denote the weight lattice (resp. the set of roots) of $E_6,E_7,E_8$ by 
$\Lambda_6,\Lambda_7,\Lambda_8$ (resp. $R_6$,$R_7$,$R_8$).
\begin{enumerate}
\item[$E_8$:] Under suitable normalization, $\Lambda_8\otimes_\Z\R$ is isometric to the Euclidean space $\R^8$. 
Let $e_1,e_2,...,e_8$ be the standard basis of $\R^8$.
The weight lattice $\Lambda_8$ is the well-known $E_8$-lattice (even coordinate system) \cite[$\mathsection8.1$]{CS88}
\begin{equation}\label{E_8}
\Lambda_8:=\{w=\sum_{i=1}^8a_ie_i\in\Z^8\cup (\Z+1/2)^8:~\sum_{i=1}^8 a_i\equiv 0~\text{mod}~2\},
\end{equation}
where $R_8$ consists of all the $240$ weights of length $\sqrt{2}$, e.g., $e_1-e_2$ and $e_2-e_3$.
This implies $R=R'$.

\item[$E_7$:]Under suitable normalization, $\Lambda_7\otimes_\Z\R$ is the orthogonal complement in $\Lambda_8\otimes_\Z\R$
 of any root $\alpha\in R_8$. The set of roots $R_7$ is just the set of $126$ elements in $R_8\cap (\Lambda_7\otimes_\Z\R)$
and $\Lambda_7$ is the orthogonal projection of $\Lambda_8$ to $\Lambda_7\otimes_\Z\R$. To see the last point,
let $\bG'$ be a complex simple Lie group of type $E_8$. 
By \cite[Ch. 8]{Ad96}, $\bG'$ contains a subgroup $\SL_2(\C)$ such that 
the identity component of its centralizer $\bG'$, denoted by $\bG$,
is of type $E_7$ and contains the center of $\SL_2(\C)$. This implies $\bG$
is simply connected. As $\bG'$ is also simply connected, the embedding $\bG\hto\bG'$ induces 
a surjective map of the weight lattices $\Lambda_8\twoheadrightarrow \Lambda_7$. This is an explicit 
description of $\Lambda_7$.
For example, we take $\alpha=e_1-e_2$. If $w=a_1e_1+\cdots+a_8e_8\in\Lambda_8$, 
then the orthogonal projection along $\alpha$ is 
\begin{equation}\label{proj7}
w':=(\frac{a_1+a_2}{2})e_1+(\frac{a_1+a_2}{2})e_2+a_3e_3+\cdots+a_8e_8.
\end{equation}
Suppose $||w'||=\sqrt{2}$ and the $a_i$ are integers, then $a_1+a_2$ must be even by
the integral equation $||w'||^2=2$ and 
we see by (\ref{proj7}),(\ref{E_8}) that $w'\in\Lambda_8$. Thus, $w'\in R_8\cap (\Lambda_7\otimes_\Z\R)=R_7$.
Suppose $||w'||=\sqrt{2}$ and all the $a_i$ are half integers, write $a_i=b_i/2$. 
We have two cases. First, if $b_1+b_2$ is even, then $w'\in\Lambda_8$ by (\ref{proj7}),(\ref{E_8}). Thus, $w'\in R_8\cap (\Lambda_7\otimes_\Z\R)=R_7$.
Second, if $b_1+b_2$ is odd, then $||w'||^2=2$ is equivalent to
\begin{equation}
(b_1+b_2)^2+2(b_3^2+\cdots+ b_8^2) = 16,
\end{equation}
which is impossible. Therefore, any weight $w'\in\Lambda_7$ of length $\sqrt{2}$ belongs to $R_7$.
We conclude that $R=R'$.

\item[$E_6$:] The strategy is similar to $E_7$. Under suitable normalization, 
$\Lambda_6\otimes_\Z\R$ is the orthogonal complement  in $\Lambda_8\otimes_\Z\R$ of a suitable chosen 
pair of roots $\alpha,\beta\in R_8$. The set of roots $R_6$ is just the set of $72$ elements in $R_8\cap (\Lambda_6\otimes_\Z\R)$
and $\Lambda_6$ is the orthogonal projection of $\Lambda_8$ to $\Lambda_6\otimes_\Z\R$.
To see the last point,
let $\bG'$ be a complex simple Lie group of type $E_8$, it contains a subgroup $\SL_3(\C)$. 
By \cite[Ch. 8]{Ad96}, the identity component of the centralizer of $\SL_3(\C)$ in $\bG'$, denoted by $\bG$,
is of type $E_6$ and contains the center of $\SL_3(\C)$. This implies $\bG$
is simply connected. As $\bG'$ is also simply connected, the embedding $\bG\hto\bG'$ induces 
a surjective map of the weight lattices $\Lambda_8\twoheadrightarrow \Lambda_6$. This is an explicit 
description of $\Lambda_6$.
For example, we take $\alpha=e_1-e_2,\beta=e_2-e_3$. If $w=a_1e_1+\cdots+a_8e_8\in\Lambda_8$, 
then the orthogonal projection along $\alpha,\beta$ is 
\begin{equation}\label{proj6}
w':=(\frac{a_1+a_2+a_3}{3})e_1+(\frac{a_1+a_2+a_3}{3})e_2+(\frac{a_1+a_2+a_3}{3})e_3+a_4e_4+\cdots+a_8e_8.
\end{equation}
Suppose $||w'||=\sqrt{2}$ and the $a_i$ are integers, then $a_1+a_2+a_3$ must be divisible by $3$ for 
the integral equation $||w'||^2=2$ and 
we see by (\ref{proj6}),(\ref{E_8}) that $w'\in\Lambda_8$. Thus, $w'\in R_8\cap (\Lambda_6\otimes_\Z\R)=R_6$.
Suppose $||w'||=\sqrt{2}$ and all the $a_i$ are half integers, write $a_i=b_i/2$. 
We have two cases. First, if $b_1+b_2+b_3$ is divisible by $3$, then $w'\in\Lambda_8$ by (\ref{proj6}),(\ref{E_8}). 
Thus, $w'\in R_8\cap (\Lambda_6\otimes_\Z\R)=R_6$.
Second, if $b_1+b_2+b_3$ is not divisible by $3$, then $||w'||^2=2$ is equivalent to
\begin{equation}
(b_1+b_2+b_3)^2+3(b_4^2+\cdots+ b_8^2) = 24,
\end{equation}
which is impossible. Therefore, any weight $w'\in\Lambda_6$ of length $\sqrt{2}$ belongs to $R_6$.
We conclude that $R=R'$.

\end{enumerate}

\subsubsection{$\mathfrak{g}=F_4$}
Under suitable normalization, $\Lambda\otimes_\Z\R$ is isometric to the Euclidean space $\R^4$. 
The weight lattice $\Lambda$ is the lattice generated by the standard basis $e_1,e_2,e_3,e_4$ and 
$\frac{e_1+e_2+e_3 +e_4}{2}$ (same as $B_4$).
The set of roots $R=\{\pm e_i,\pm e_i\pm e_j,\frac{\pm e_1 \pm e_2\pm e_3\pm e_4}{2}:~1\leq i\neq j\leq n\}$.

Since  $R$ consists of all weights of length $1$ or $\sqrt{2}$, we obtain $R=R'$.

\subsubsection{$\mathfrak{g}=G_2$}
The weight lattice $\Lambda$ is generated by the set of roots $R$.
By observing the root system of $G_2$, one sees that 
the shortest nonzero weights of $\Lambda$ belong to $R$. 
Hence, the short roots $R'^\circ$ of $R'$ belong to $R$. 
Since $R'^\circ$ determines $R'$, we obtain $R=R'$.

\subsubsection{Theorem \ref{thm6}(ii)}
It remains to prove Theorem \ref{thm6}(ii).
Let $\bG'$ be a connected complex simple Lie group with Lie algebra $\mathfrak{g}'\in\{E_7,E_8,F_4\}$. 
Embed $\bG'$ in some
$\GL_N(\C)$. When $\bG'$ is of type $E_7$, it contains an equal rank subgroup $\bG$ of type $A_7$;
when $\bG'$ is of type $E_8$, it contains an equal rank subgroup $\bG$ of type $A_8$ (resp. $D_8$);
when $\bG'$ is of type $F_4$, it contains an equal rank subgroup $\bG$ of type $B_4$ \cite[Table 5]{GOV94}.
Then we obtain an embedding $f:\bG\hto\bG'\hto\GL_N(\C)$.

Pick a maximal torus $\bT$ of $\bG$. It is also a maximal torus of $\bG'$.
If the normalizer $\bN:=N_{\GL_N(\C)}(\bT)$ preserves the set of roots $R$ of $\bG$, then 
the Weyl group $W'$ of $\bG'$ preserves $R$.
Since $W'$ is orthogonal on $\Lambda\otimes_\Z\R$, $W'$ acts on the set $P$ of 
\emph{Weyl chambers} of $R\subset\Lambda\otimes_\Z\R$ \cite[$\mathsection14$]{FH91}. 
Each Weyl chamber corresponds to a set of positive roots $R^+$ of $R$. 
The order of $P$ is equal to the order of 
$W$, the Weyl group of $\bG$. 

Let $R^+\in P$ and $H$ be the stabilizer of $R^+$ in $W'$.   
Since $W'$ is faithful and orthogonal on $\Lambda\otimes_\Z\R$, $H$ is an automorphism group 
of the \emph{Dynkin diagram} \cite[Table 1]{GOV94} of $\bG$ (or $\mathfrak{g}$).
As $W$ is simply transitive on $P$, we obtain $|H|=|W'|/|W|$. \\

\begin{minipage}{0.3\textwidth}  
\begin{center}
\begin{tabular}{|c|c|} \hline
$\mathfrak{g}'$ & $|W'|$\\ \hline
$E_7$ & 2903040 \\ \hline
$E_8$ & 696729600 \\ \hline
$F_4$ &  1152 \\ \hline
\end{tabular}
\end{center}
\end{minipage}
\begin{minipage}{0.3\textwidth}  
\begin{center}
\begin{tabular}{|c|c|c|} \hline
$\mathfrak{g}$ & $|W|$ & Dynkin diagram\\ \hline
$A_7$ & $8!=40320$ & 
\begin{tikzpicture}[scale=.25]
    \draw (-1,0) ;
    \foreach \x in {0,...,6}
    \draw[xshift=\x cm,thick] (\x cm,0) circle (.3cm);
    \draw[thick] (0.3 cm,0) -- +(1.4 cm,0);
    \foreach \y in {1.15,...,5.15}
    \draw[xshift=\y cm,thick] (\y cm,0) -- +(1.4 cm,0);
  \end{tikzpicture} \\ \hline
	
$A_8$ & $9!=362880$ & 
\begin{tikzpicture}[scale=.25]
    \draw (-1,0) ;
    \foreach \x in {0,...,7}
    \draw[xshift=\x cm,thick] (\x cm,0) circle (.3cm);
    \draw[thick] (0.3 cm,0) -- +(1.4 cm,0);
    \foreach \y in {1.15,...,6.15}
    \draw[xshift=\y cm,thick] (\y cm,0) -- +(1.4 cm,0);
  \end{tikzpicture}\\ \hline
	
$B_4$ &  $4!\times 2^4=384$ & 
\begin{tikzpicture}[scale=.25]
    \draw (-1,0) ;
    \foreach \x in {0,...,2}
    \draw[xshift=\x cm,thick] (\x cm,0) circle (.3cm);
    \draw[xshift=5 cm,thick] (1 cm, 0) circle (.3 cm);
    \draw[thick] (0.3 cm,0) -- +(1.4 cm,0);
    \foreach \y in {0.15,...,1.15}
    \draw[xshift=\y cm,thick] (\y cm,0) -- +(1.4 cm,0);
    \draw[thick] (4.3 cm, .1 cm) -- +(1.4 cm,0);
    \draw[thick] (4.3 cm, -.1 cm) -- +(1.4 cm,0);
		\draw
(5.1 cm,0) --++ (120:0.4)
(5.1 cm,0) --++ (-120:0.4);
  \end{tikzpicture}\\ \hline
	
$D_8$ &  $8!\times 2^7=5160960$ & 
 \begin{tikzpicture}[scale=.25]
    \draw (-1,0);
    \foreach \x in {0,...,5}
    \draw[xshift=\x cm,thick] (\x cm,0) circle (.3cm);
    \draw[xshift=10 cm,thick] (30: 17 mm) circle (.3cm);
    \draw[xshift=10 cm,thick] (-30: 17 mm) circle (.3cm);
    \draw[thick] (0.3 cm,0) -- +(1.4 cm,0);
    \foreach \y in {1.15,...,4.15}
    \draw[xshift=\y cm,thick] (\y cm,0) -- +(1.4 cm,0);
    \draw[xshift=10 cm,thick] (30: 3 mm) -- (30: 14 mm);
    \draw[xshift=10 cm,thick] (-30: 3 mm) -- (-30: 14 mm);
  \end{tikzpicture}  \\ \hline
\end{tabular}
\end{center}
\end{minipage}
~\\

For any embedding $f:\bG\hto\GL_N(\C)$ in the first paragraph, we have $|H|>2$ by the table above.
This contradicts that the automorphism group of the Dynkin diagram of $A_7,A_8,B_4,D_8$ has order 
at most $2$. We conclude that $\bN$ does not preserve the set of roots $R$.

\begin{remark}
The exceptional complex simple Lie group $\bG'$ of type $G_2$ contains $\SL_3(\C)$ (type $A_2$) as an equal rank subgroup.
The arguments above do not contradict Theorem \ref{thm6}(i) for $\bG=\SL_3(\C)$
because the orders of the Weyl groups of $G_2$ and $A_2$ are 
respectively $12$ and $6$ and $\frac{12}{6}=2$ is the order of the automorphism of the Dynkin diagram of $A_2$.
\end{remark}

\section*{Acknowledgments} I would like to thank Gabor Wiese for his interests and comments on the paper.
I would like to thank the referee for his/her constructive comments.

\vspace{.1in}
\noindent Yau Mathematical Sciences Center, Tsinghua University, Haidian District, Beijing 100084, China\\
Emails: \url{pslnfq@gmail.com}, \url{pslnfq@tsinghua.edu.cn}

\end{document}